\mag=\magstep1
\input epsf
\documentclass[10pt,a4paper]{amsart}
\usepackage{amssymb,latexsym}

\def\Cal{\mathcal}
\def\subsetne{\underset{\neq}\subset}

\def\len{\operatorname{len}}
\def\op{\operatorname{op}}
\def\ot{\leftarrow}

\def\<<{\langle}
\def\>>{\rangle}

\numberwithin{equation}{section}
\newtheorem{theorem}{Theorem}[section]
\newtheorem{proposition}[theorem]{Proposition}
\newtheorem{corollary}[theorem]{Corollary}
\newtheorem{definition}[theorem]{Definition}
\newtheorem{conjecture}[theorem]{Conjecture}
\newtheorem{remark}[theorem]{Remark}
\newtheorem{lemma}[theorem]{Lemma}

\begin{document}

\title{
Mixed Tate motives and multiple zeta values.
}

\author{Tomohide Terasoma}
\address{Department of Mathematical Science, University of Tokyo, Komaba 3-8-1, 
Meguro, Tokyo 153 , Japan}
\email{terasoma@ms.u-tokyo.ac.jp}


\date{\today }
%

\maketitle

\makeatletter
\renewcommand{\@evenhead}{\tiny \thepage \hfill  T.TERASOMA \hfill}
\renewcommand{\@oddhead}{\tiny \hfill  FUNDAMENTAL MOTIVES AND MULTIPLE
ZETA VALUES 
 \hfill \thepage}
\makeatother

\section{Introduction}
\label{sec:introduction}

Let $l$ be a positive integer and $k_1, \dots , k_l$
be integers such that $k_i \geq 1$ for $i=1, \dots ,l-1$ and
$k_l \geq 2$. We define the multiple zeta value $\zeta (k_1, \dots ,k_l)$
of index $(k_1, \dots ,k_l)$ as
$$
\zeta (k_1, \dots ,k_l)=\sum_{m_1 < \cdots < m_l}\frac{1}
{m_1^{k_1}\cdots m_l^{k_l}}.
$$
Set $h_0=0$ and $h_i = \sum_{j=1}^ik_i$. 
For an index $(k_1, \dots, k_l)$, the number $n = \sum_{i=1}^lk_i$
is called the weight of the index. The weight of
the multiple zeta value $\zeta (k_1, \dots ,k_l)$ is defined to be $n$.
We have
the following integral expression of $\zeta (k_1, \dots ,k_l)$:
$$
\zeta (k_1, \dots ,k_l)=
(-1)^l
\int_{\Delta_n}\prod_{i=1}^l
[\frac{dx_{h_{i-1}+1}}{x_{h_{i-1}+1}-1}\wedge
\prod_{j=1}^{k_{i}-1}\frac{dx_{h_{i-1}+j+1}}{x_{h_{i-1}+j+1}}]
$$
Here $\Delta_n$ is the topological cycle defined by
$\{0< x_1 < \dots < x_n < 1\}$.
We define the subspace $L_n$ of $\bold R$ by
$$
L_n= \sum_{\substack{k_1+\cdots +k_l=n \\ k_i \geq 1, k_l \geq 2}}
\zeta (k_1, \dots ,k_l)\bold Q.
$$
By the shuffle relation, we have $L_i \cdot L_j \subset L_{i+l}$
and the space $\Cal A = \oplus_{n=0}^\infty L_{n}$
is a graded algebra. Zagier \cite{Z} conjectured 
the following dimension formula. (See also \cite{BK}.)
For the motivic interpretation, see \cite{D}.
\begin{conjecture}[Zagier]
$$
\dim L_n = d_n,
$$
where $d_i$ is defined by the following inductive formula.
$$
d_0=1, d_1=0, d_2=1, d_{i+3}=d_{i+1}+d_i \text{ for } i \geq 0
$$
\end{conjecture}
The generating function of $d_i$ is given by 
$$
\sum_{i=0}^\infty d_it^i =\frac{1}{1-t^2-t^3}.
$$
The main theorem of this paper is the following.
\begin{theorem}[Main Theorem]
$\dim L_n \leq d_n$.
\end{theorem}
The same theorem is announced in \cite{G} and \cite{G2}.
Let us explain the outline of the proof. We construct 
a pair of varieties $(Y^0, \bold B^0)$ such that
the periods of the relative cohomology 
$H^n = H^n(Y^0,j_! \bold Q_{Y^0 - \bold B^0})$
are expressed
as $\bold Q$-linear combinations of multiple zeta values.
Here $j$ is the natural inclusion $Y^0 - \bold B^0 \hookrightarrow Y^0$.
We show that the extension class of
\begin{align}
\label{eqn:length one}
0 \to W_{k+1}(H^n)/W_{k+2}(H^n) & \to W_{k}(H^n)/W_{k+2}(H^n) \\
& \to W_{k}(H^n)/W_{k+1}(H^n) \to 0 \nonumber
\end{align}
vanishes in the category of mixed Hodge structures for all $k$.
We use the the following equality for the extension of 
mixed motives after Marc Levine \cite{L}.
\begin{equation}
\label{eqn:mixed ext}
Hom_{DTM(\bold Q)}(\bold Q, \bold Q(i)[1])= K_{2i-1}(\bold Q) \otimes \bold Q
\end{equation}
In particular, we have 
$Hom_{DTM(\bold Q)}(\bold Q, \bold Q(1)[1])=\bold Q^{\times} \otimes \bold Q$.
By using the compatibility in Proposition \ref{prop:compatiblity ext one},
the extension class (\ref{eqn:length one}) is known to 
split in the category of mixed motives.
Using this fact and equality (\ref{eqn:mixed ext}), 
we prove that $H^n(Y^0,j_! \bold Q_{Y^0 - \bold B^0})$
is a subquotient of a direct sum of typical objects 
in the category of mixed motives.

Let us explain the contents of this paper. In Section \ref{sec:var ass fund gp},
we construct by a succession of blowing ups
a variety $Y^0$ and its closed subvariety $\bold B^0$
associated to multiple zeta values.
We compute the relative de Rham cohomology 
$H^i(j_! K_{Y^0-\bold B^0,DR})$. 
Since $H^i(j_! K_{Y^0-\bold B^0,DR})=0$
for $i \neq n$, the cone 
$Cone(\bold Q_{Y^0} \to \bold Q_{\bold B^0})[n]$
defines an abelian object $\Cal A_{TM}$ in $DTM_{\bold Q}$.
In Section \ref{sec:generator MTM}, 
we recall the definition of the abelian category of
mixed Tate motives and prove several propositions concerning
generators of abelian sub-categories in $\Cal A_{TM}$.
We also prove the compatibility of
the cycle map for mixed Tate motives and 
the map $ch$ for the extensions of mixed Tate Hodge structures.
In Section \ref{sec:periods MTS MZV}, 
we compute the periods of the relative cohomologies
defined in Section \ref{sec:var ass fund gp}.
We claim that the periods of the relative cohomologies are expressed
by multiple zeta values.
In the last section, we prove the Main Theorem.

The author would like to express his thanks to M.Matsumoto, T.Saito, M.Kaneko
and B. Kahn for discussions. He also would like to express his thanks
to M.Levine for letting him know the references about 
the category of mixed Tate motives.

\section{Varieties associated to multiple zeta values}
\label{sec:var ass fund gp}
\subsection{Successive blowing up of affine spaces}
\label{subsec:successive}
In this section, we introduce a successive blowing up of 
$(\bold A^1)^n$.

Let $(\bold A^1)^n = \{(x_1, \dots , x_n) \}$ be a product of $\bold A^1$.
We define a sequence of subvarieties $Z_i, W_i$ for $i=0, \dots , n-2$
by
\begin{align*}
Z_i &  = \{ (x_1, \dots , x_n) \mid x_1 = \cdots = x_{n-i} = 0\}, \\
W_i & = \{ (x_1, \dots , x_n) \mid x_{i+1} = \cdots = x_n = 1\}, 
\end{align*}
We define sequences of blowing ups and their centers

\begin{center}
\vskip 0.1in
\begin{tabular}{ccccccccc}
$(\bold A^1)^n = X_0$ & $\ot$ &  $X_1$ &  $\ot$ & $\cdots$ 
& $\ot$ & $X_{n-2}$ & $\ot$ & $X_{n-1}=X$ \\
$\cup$ & & $\cup$ & & & & $\cup$ \\
$Z_0^{pr}$  & & $Z_1^{pr}$ & & & &   $Z_{n-2}^{pr}$ \\
\end{tabular}
\vskip 0.1in
\begin{tabular}{ccccccccc}
$X = Y_0$ & $\ot$ &  $Y_1$ & $\ot$ & $\cdots$ & $\ot$ 
& $Y_{n-2}$ & $\ot$ & $Y_{n-1}=Y$ \\
$\cup$ & & $\cup$ & & & & $\cup$ \\
$W_0^{pr}$ & & $W_1^{pr}$ & & & & $W_{n-2}^{pr}$ \\
\end{tabular}
\vskip 0.1in
\end{center}

\noindent
by the following procedure.
\begin{enumerate}
\item
$X_0$ (resp. $Y_0$) is equal to $(\bold A^1)^n$ (resp. $X$). 
\item
Let $Z_i^{pr}$ ($W_i^{pr}$) be the proper transform of $Z_i$ (resp. $W_i$)
under the morphism $X_i \to X_0$ (resp, $Y_i \to X_0$).
\item 
Let $X_{i+1}$ (resp. $Y_{i+1}$) be the blowing up of $X_i$
(resp. $Y_i$) with the center $Z_i^{pr}$ (resp. $W_i^{pr}$)
for $i=1, \dots ,n-2$.
\end{enumerate}

We define divisors $D_i^{\epsilon}$ ($\epsilon = 0, 1, 1 \leq i \leq n$)
of $(\bold A^1)^n$ by
$$
D_i^{\epsilon} =\{ x_i = \epsilon \}
$$
and
$$
D= D_2^{0}\cup \cdots  \cup D_n^0 \cup 
D_1^{1}\cup \cdots  \cup D_{n-1}^1  
$$
and the proper transform of $D$ under the morphism $X_i \to X_0$
(resp. $Y_i \to X_0$, $X\to X_0$, $Y \to X_0$) is denoted by $D_{X_i}$ 
(resp. $D_{Y_i}$, $D_{X}$, $D_{Y}$).
The proper transform of $D_j^{\epsilon}$ under the morphism
$X_i \to X_0$
(resp. $Y_i \to X_0$, $X \to X_0$, $Y \to X_0$) 
is denoted by $D_{j,X_i}^{\epsilon}$ 
(resp. $D_{j, Y_i}^{\epsilon}$, $D_{j, X}^{\epsilon}$, $D_{j, Y}^{\epsilon}$).
The open subvariety $X_i - D_{X_i}$ (resp. $Y_i - D_{Y_i}$, $X-D_X$, $Y-D_Y$)
is denoted by $X_i^0$ (resp. $Y_i^0$, $X^0$, $Y^0$).

We define $Z_i^0$ and $W_i^0$ by $Z_i^{pr}\cap X_i^0$ and $W_i^{pr}\cap Y_i^0$,
respectively.
It is easy to see that the divisors $D_{1, X_i}^0 , \dots ,D_{n-i, X_i}^0$
(resp. $D_{i+1, Y_i}^1 , \dots ,D_{n, Y_i}^1$)
are normal crossing and
$Z_{i}^{pr} = \cap_{j=1}^{n-i}D_{j, X_i}^0$
(resp. $W_{i}^{pr} = \cap_{j=i+1}^{n}D_{j, Y_i}^1$). Therefore the inverse
image of $X_i^0$ (resp. $Y_i^0$) 
under the morphism $X_{i+1} \to X_i$ (resp. $Y_{i+1} \to Y_i$) is isomorphic
to $X_i^0$ (resp. $Y_i^0$). 
We identify $X_{i}^0$ (resp. $Y_i^0$) as an open set of 
$X_{i+1}$ (resp. $Y_{i+1}$)
via this isomorphism. Let $E_i^0 = X_{i+1}^0 - X_i^0$ 
(resp. $F_i^0 = Y_{i+1}^0 - Y_i^0$).

In general, for a system of coordinate $(y_1, \dots , y_k)$ of 
$(\bold A^1)^k$, we consider the same procedure as before and
the resultant variety $X$ and $ Y$ are denoted by $X(y_1, \dots ,y_k)$
and $Y(y_1, \dots , y_k)$, respectively.
The divisors corresponding to
$D_Y$, $D_{i,X}^{\epsilon}$ and $D_{i,Y}^{\epsilon}$ are denoted by
$D_Y(y_1, \dots ,  y_k)$, $D_{y_i,X}^{\epsilon}(y_1, \dots , y_k)$ and 
$D_{y_i,Y}^{\epsilon}(y_1, \dots , y_k)$, respectively.
The following propositions can be easily checked.

\begin{proposition}
\begin{enumerate}
\item
The divisors $D_{j,X_i}^0$ ($j=n-i+2, \dots , n$)
and $D_{j,X_i}^1$ ($j=n-i+1, \dots , n$)
transversally intersect with $Z_i^{pr}$.
The divisors $D_{j,Y_i}^1$ ($j=1, \dots , i-1$)
and $D_{j,Y_i}^0$ ($j=1, \dots , i$)
transversally intersect with $W_i^{pr}$.
\item
The divisor $D_{j,X_i}^1$ does not intersect with $Z_i^{pr}$
for $j=1, \dots , n-i$.
The divisor $D_{n-i+1,X_i}^0$ does not intersect with $Z_i^{pr}$.
\item
The divisor $D_{j,Y_i}^0$ does not intersect with $W_i^{pr}$
for $j=i+1, \dots , n$.
The divisor $D_{i,Y_i}^1$ does not intersect with $W_i^{pr}$.
\end{enumerate}
\end{proposition}
\begin{proposition}
\label{prop:induct z}
The variety $Z_i^{pr}$ (resp. $W_i^{pr}$) is naturally identified with
$X(x_{n-i+1}, \dots , x_n)$ (resp. $Y(x_1, \dots , x_i)$).
\begin{enumerate}
\item
Under the identification $Z_i^{pr}=X(x_{n-i+1}, \dots ,x_n)$,
the intersection 
\linebreak
$D_{x_j,X_i}^0 \cap Z_i^{pr}$
(resp. $D_{x_j,X_i}^1 \cap Z_i^{pr}$) is identified with
$D_{x_j, X}^0(x_{n-i+1}, \dots , x_n)$ for $n-i+2 \leq j \leq n$
(resp. $D_{x_j,X}^1(x_{n-i+1}, \dots , x_n)$ for $n-i+1 \leq j \leq n-1$).
\item
Under the identification $W_i^{pr}=Y(x_1, \dots ,x_i)$, 
the intersection $D_{x_j,Y_i}^0 \cap W_i^{pr}$
(resp. $D_{x_j,Y_i}^1 \cap W_i^{pr}$) is identified with
$D_{x_j, Y}^0(x_1, \dots , x_i)$ for $2 \leq j \leq i$
(resp. 
$D_{x_i,Y}^1(x_1, \dots , x_i)$ for $1 \leq j \leq i-1$).
\end{enumerate}
\end{proposition}
\begin{corollary}
\label{cor:z open}
\begin{align*}
 Z_i^0 = & Z_i^{pr} -\cup_{j=n-i+2}^nD_{j,X_i}^0(x_{n-i+1}, \dots ,x_n) \\
& -\cup_{j=n-i+1}^{n-1}D_{j,X_i}^1(x_{n-i+1}, \dots ,x_n),
\\
W_i^0 = & W_i^{pr} -\cup_{j=2}^iD_{j,Y_i}^0(x_{1}, \dots ,x_i)
-\cup_{j=1}^{i-1}D_{j,Y_i}^1(x_{1}, \dots ,x_i).
\end{align*}
\end{corollary}

Let $E_i$ (resp. $F_i$) be the exceptional divisor of the blowing up
$X_{i+1} \to X_i$ (resp. $Y_{i+1} \to Y_i$). 
The the morphism $\pi_i:E_i \to Z_i^{pr}$ ($\tau_i:F_i \to W_i^{pr}$)
is a $\bold P^{n-i-1}$-bundle and the intersections
$D_{j,X_{i+1}}^0 \cap E_i$  ($1 \leq j \leq n-i$)(resp. 
$D_{j,Y_{i+1}}^1 \cap F_i$ ($i+1 \leq j \leq n$))
are horizontal 
families of independent hyperplanes for the morphism 
$E_i \to Z_i^{pr}$ (resp. $F_i \to W_i^{pr}$).
Therefore $E_i-\cup_{j=2}^{n-i} D_{j,X_{i+1}}^0$ 
(resp. $F_i-\cup_{j=i+1}^{n-1} D_{j,X_{i+1}}^1$)
is a $\bold A^1 \times
(\bold G_m)^{n-2-i}$-bundle over $Z_i^{pr}$ (resp. $W_i^{pr}$). 
By Corollary \ref{cor:z open},
the morphism $\pi_i$ (resp. $\tau_i$)
induces a morphism $E_i^0 \to Z_i^0$ (resp. $F_i^0 \to W_i^0$), which is
also denoted by $\pi_i$ (resp. $\tau_i$).

We introduce an open set $U_i$ of $Z_i^0$ and its neighborhood $N_i$.
We use them and their blowing ups to compute the de Rham cohomology
of $X^0$ in Proposition \ref{prop:de Rham cohom blow up}.

By Proposition \ref{prop:induct z} and Corollary \ref{cor:z open}, 
$Z_i^0$ contains an openset $U_i$ defined by
\begin{align*}
U_i = \{ (x_{n-i+1}, \dots , x_n) \mid
& x_k \neq 0 (\text{ for } n-i+1 \leq k \leq n ), \\
& x_k \neq 1 (\text{ for } n-i+1 \leq k \leq n-1) 
\}. 
\end{align*}

We give a description of the restriction of $\pi_i$ to $U_i$. Let 
$$
N_i =\{(x_1, \dots , x_n) \mid (x_{n-i+1}, \dots ,x_n) \in U_i \}.
$$
The variety $U_i$ can be identified with the closed subvariety
$\{(x_1, \dots ,x_n)\in N_i \mid x_1= \cdots = x_{n-i} =0\}$ of $N_i$.
Since $N_i$ does not intersect with $Z_0, \dots ,Z_{i-1}$ in
$X_0$, $N_i$ can be identified with a subvariety in $X_i$
via the blowing up procedure. Moreover we have
$N_i \cap Z_i^0 = U_i$ in $X_i$. Let $Bl_{U_i}(N_i)$ be
the blowing up of $N_i$ along the center $U_i$. Then 
$Bl_{U_i}(N_i)$ can be identified with an open set of
$X_{i+1}$, and we have
the following cartesian diagram.

\begin{center}
\vskip 0.1in 
\begin{tabular}{ccc}
$X_{i+1}$ & $\to$ & $X_i$ \\
$\cup$ &  & $\cup$ \\
$Bl_{U_i}(N_i)$ & $\to$ & $N_i$ \\
\end{tabular}
\vskip 0.1in
\end{center}

Let $Bl_{U_i}(N_i)^0 = Bl_{U_i}(N_i)-\cup_{j=2}^{n-i}D_{j,X_{i+1}}^0
-\cup_{j=1}^{n-i}D_{j,X_{i+1}}^1$. Then $Bl_{U_i}(N_i)^0=
X_{i+1}^0 \cap Bl_{U_i}(N_i)$. We put 
$E_{i,U_i}^0=Bl_{U_i}(N_i)^0 \cap E_i^0$.
Then we have the following commutative diagram.

\vskip 0.1in
\begin{center}
\begin{tabular}{ccc}
$E_i^0$ & $\to$ & $Z_i^0$ \\
$\uparrow$ & & $\uparrow$ \\
$E_{i,U_i}^0$ & $\to$ & $U_i$ \\  
\end{tabular}
\end{center}
\vskip 0.1in

\begin{proposition}
\label{prop:excep div}
We introduce coordinates $\xi_i$ with $\xi_i x_2 = x_i$ for
$i=1,3,4,$ $\dots ,n-i$.
\begin{enumerate}
\item
The variety $Bl_{U_i}(N_i)^0$ is isomorphic to 
\begin{align*}
\{  (\xi_1, & \xi_3, \dots ,\xi_{n-i}, x_2, x_{n-i+1}, \dots ,x_n) \mid
 \xi_i \neq 0 (\text{ for }3 \leq i \leq n-i), \\
& x_2\xi_i \neq 1 (\text{ for }3 \leq i \leq n-i \text{ or } i=1), 
 x_i \neq 0 (\text{ for }n-i+1 \leq i \leq n), \\
& x_i \neq 1 (\text{ for }n-i+1 \leq i \leq n-1 \text{ or } i=2) \}.
\end{align*}
\item
The variety $E_{i,U_i}^0$ is isomorphic to 
the subvariety of $Bl_{U_i}(N_i)^0$ defined by $x_2=0$
under the coordinate given as above. More explicitly,
it is isomorphic to
\begin{align*}
\{(\xi_1, \xi_3,& \dots ,\xi_{n-i},x_{n-i+1}, \dots ,x_n) \mid
 \xi_i \neq 0 (\text{ for }3 \leq i \leq n-i), \\
& x_i \neq 0 (\text{ for }n-i+1 \leq i \leq n), 
 x_i \neq 1 (\text{ for }n-i+1 \leq i \leq n-1 ) \}. 
\end{align*}
\end{enumerate}
\end{proposition}
We also use the following proposition for the computation of de Rham
cohomology in the next subsection.
\begin{proposition}
\label{prop:bold U}
The open set $Bl_{U_i}(N_i)^0 -E_{i,U_i}^0$ is isomorphic to
\begin{align*}
\bold U = \{(x_1, \dots , x_n) \mid
& x_i \neq 0 (\text{ for }2 \leq i \leq n), \\
& x_i \neq 1 (\text{ for }1 \leq i \leq n-1) \}. 
\end{align*}
\end{proposition}

\subsection{Computation of de Rham cohomology}
In this section, we compute the de Rham cohomology of $X^0$ and $Y^0$.
Let $K$ be the field of definition $\bold Q$.

We define a set $\Cal S$ of the 
differential forms on $X_0^0$ by
$$
\Cal S = \{ \frac{dx_i}{x_i} \ (i= 2, \dots , n),
 \frac{dx_i}{x_i-1} \ (i= 1, \dots , n-1) \}
$$
Let $1 \leq i_1 < \cdots < i_k \leq n$, 
$\displaystyle\frac{dx_{i_r}}{x_{i_r}-e_r}$
be an element in $\Cal S$ for $1 \leq r \leq k$. A differential form
\begin{equation}
\label{eqn:diff notation}
\omega = \omega (i_1,\dots ,i_k; e_1, \dots , e_k)=
\frac{dx_{i_1}}{x_{i_1}-e_1}\wedge \cdots \wedge
\frac{dx_{i_k}}{x_{i_k}-e_k}
\end{equation}
is said to begin with type 1 (resp. end by type 0) if 
$e_1=1$ (resp. $e_k=0$).
The vector space of $k$-forms generated by 
$\omega (i_1, \dots ,i_k ;e_1, \dots , e_k)$ with 
$e_1 = 1$ (resp. $e_1=1$ and  $e_k=0$)
is denoted by $V^k_{1}$ (resp. $V^k_{10}$).
If $k=0$, we define $V_1^0= V_{10}^0=K$.
The main theorem of this section is as follows.

\begin{theorem}
\label{thm:de Rham cohom bl up}
\begin{enumerate}
\item
$
H^k_{DR}(X^0) \simeq V^k_{1}
$.
\item
$
H^k_{DR}(Y^0) \simeq V^k_{10}
$.
\end{enumerate}
\end{theorem}

Let $V_{k,p}$ be the subspace of $k$-forms $H^k_{DR}(\bold U)$
of $\bold U$ generated by
\linebreak
$\omega (i_1, \dots, i_k;e_1, \dots , e_k)$, where
(1) there exists $r$ such that $i_r=p$ and $e_r=1$, and
(2) $e_q=0$ for $q< r$.
For a formal sum $a = \sum_{J\subset [1,n]} a_J J$ 
($a_J \in \bold Q$) of 
simplices $J$ in $[1,n]$, we define
$\omega^0 (a,x)$ by
\begin{equation}
\label{eqn:mult simplex diff}
\omega^0 (a,x) = \sum_J a_J \prod_{j \in J}\frac{dx_j}{x_j}.
\end{equation}
The boundary operator is written as $\partial$ and
the $k$-th face operator is denoted by ``$-\{ k\}$''.
For $2\leq q < p \leq n+1$,
we define a subspace $V_{k,p,q}$ by the subspace of
$k$-form generated by 
$\omega^0 (\partial J,x)\wedge \frac{dx_p}{x_p-1} \wedge \eta$
(resp. $\omega^0 (\partial J,x)$),
if $p \leq n-1$ (resp. $p=n$ or $p=n+1$)
where $J \subset [2, q]$ and $\eta$ a product of
$\displaystyle\frac{dx_j}{x_j-e_j} \in \Cal S$ with $j > p$.
Since 
$V_{k,p,2}=V_{k,p}$ and
$V_1^k = V_{k,1}\oplus V_{k,2} \oplus \oplus_{p=3}^{n-1}V_{k,p,2}$,
the following proposition implies Theorem \ref{thm:de Rham cohom bl up}.1.

\begin{proposition}
\label{prop:de Rham cohom blow up}
\begin{enumerate}
\item
The natural map
$H^k(X_i^0 ) \to H^k(\bold  U )$ is injective.
\item
Under the above injection, $H^k(X_i^0 )$ ($i=1, \dots, n-1$)
is identified with
$$
\oplus_{p=1}^{n-i+1}V_{k,p}\oplus \oplus_{p=n-i+2}^{n+1} V_{k,p,n-i+1}.
$$
\end{enumerate}
\end{proposition}
\begin{proof}
We prove the proposition by induction of $\dim X$ and $i$.
For $i=1$, we can prove the proposition directly.
We consider the following commutative diagram whose rows are exact:

\begin{center}
\vskip 0.1in
\begin{tabular}{ccccc}
$H^k(X_{i+1}^0 )$ & $\to$ &
$H^k(X_i^0 )$ & $\overset{r}\to$ & $H^{k-1}(E_{i}^0)$ \\
$\downarrow$ & & $\alpha\downarrow$ & & $\downarrow\beta$ \\
$H^k(Bl_{U_i}(N_i)^0 )$ & $\to$ &
$H^k(\bold  U )$ & $\overset{r}\to$ & $H^{k-1}(E_{i,U_i}^0)$ 
\end{tabular}
\vskip 0.1in
\end{center}

By the hypothesis of the induction, $\alpha$ is injective.
By using the Leray spectral sequence and inductive hypothesis, $\beta$ is
also injective. By Proposition \ref{prop:bold U} and 
Proposition \ref{prop:excep div}
$H^k(\bold U)$, $H^k(E_{i,U_i}^0)$ and
$H^k(E_{i}^0)$ are generated by $\Cal B_{\bold U}$,
$\Cal B_{E_{i,U_i}^0}$ and $\Cal B_{E_{i}^0}$, where
\begin{align*}
& \Cal B_{\bold U} =  \{ \omega (i_1, \dots ,i_k;e_1, \dots ,e_k) 
\mid \frac{dx_{i_k}}{x_{i_k}-e_i}
\in \Cal S \}, \\
& \Cal S^{(i)}= \{ \frac{dx_i}{x_i} \ (i=n-i+1, \dots , n),
\frac{dx_i}{x_i-1} \ (i=n-i+1, \dots , n-1) \}, \\
& \Cal B_{E_{i,U_i}^0}= 
 \{ \omega^0 (J,\xi) \wedge\omega (j_1, \dots ,j_l; e_1, \dots ,e_l)\mid
J \subset [3, n-i], \frac{dx_{j_k}}{x_{j_k}-e_i} \in \Cal S^{(i)} 
\}, \\
& \Cal B_{E_{i}^0}=  \{ 
\omega^0 (J,\xi) \wedge\omega (j_1, \dots ,j_l; e_1, \dots ,e_l) \in 
\Cal B_{E_{i,U_i}^0}
\mid e_1 = 1 \}. 
\end{align*}
Here we used notations (\ref{eqn:diff notation}) and
(\ref{eqn:mult simplex diff}).
The horizontal homomorphisms $r$ are obtained by Poincare residue
with respect to the divisor $x_2=0$.
We prove that $r : H^k(X_i^0) \to H^k(E_i^0)$ 
is surjective and the kernel of $r$ is equal to
$\oplus_{p=1}^{n-i}V_{k,p }\oplus \oplus_{p=n-i+1}^{n+1} V_{k, p, n-i}$.
Let $T =(i_1, \dots ,i_l;e_1, \dots ,e_l)$ with $i_1, \dots , i_l \in [p+1,n]$.
We use the notation $\omega (T)$ in (\ref{eqn:diff notation}).
We define  $V_{k,p}^{T}$ and $V_{k,p,q}^{T}$ as the subspace of
$V_{k,p}$ and $V_{k,p,q}$ generated by
$\displaystyle\omega^0 (J,x) \wedge\frac{dx_p}{x_p-1} \wedge \omega (T)$
and 
$\displaystyle\omega^0 (\partial J,x) \wedge\frac{dx_p}{x_p-1} \wedge \omega (T)$ with
$J \in [2,q]$, respectively.
Then we have $V_{k,p} =\oplus_T V_{k,p}^{T}$ and
$V_{k,p,q} =\oplus_T V_{k,p,q}^{T}$ for $p > n-i$.
We define $V_{k,p,n-i}^{T,\xi}$ by
the subspace of $H^k(E_{i,U_i}^0)$ generated by
$$
\{\omega^0 (J,\xi) \wedge \frac{dx_p}{x_p-1} \wedge \omega (T)
\mid J \subset [3, n-i] \}.
$$
We compute the residue $r$ with respect to $x_2=0$ by the
relation $\displaystyle\frac{dx_k}{x_k} = 
\frac{d\xi_k}{\xi_k} + \frac{dx_2}{x_2}$
for $k=3, \dots , n-i$.
It is easy to see that $r(V_{k,p,n-i+1}^{T}) \subset V_{k,p,n-i}^{T,\xi}$
for $p > n-i+1$.
For $J \subset [2, n-i]$, we have
\begin{equation}
\label{eqn:residue}
r(\omega^0 (J,x)\wedge \frac{dx_p}{x_p-1} \wedge \omega (T))=
\begin{cases}
\omega^0 (\partial J, \xi )\wedge \frac{dx_p}{x_p-1} \wedge \omega (T)
& \text{ if }2 \notin J   \\
\omega^0 (J-\{ 2 \}, \xi )\wedge \frac{dx_p}{x_p-1} \wedge \omega (T)
& \text{ if }2 \in J   
\end{cases}
\end{equation}
and $r(V_1)= \cdots = r(V_{n-i})=0$.
Therefore it is enough to prove that 
$r:V_{k,n-i+1}^{T} \to V_{k,n-i+1,n-i}^{T,\xi}$ and
$r:V_{k,p,n-i+1}^{T} \to V_{k,p,n-i}^{T,\xi}$ ($p > n-i+1$) are surjective and
\begin{align*}
& ker(r:V_{k,n-i+1}^{T} \to V_{k,n-i+1,n-i}^{T,\xi})=V_{k,n-i+1,n-i}^{T}, \\
& ker(r:V_{k,p,n-i+1}^{T} \to V_{k,p,n-i}^{T,\xi})=V_{k,p,n-i}^{T}.
\end{align*}

(1)
{\bf The morphism} 
$r:V_{k,n-i+1}^{T} \to V_{k,n-i+1,n-i}^{T,\xi}$.
For $J \subset [3,n-i]$,
$\displaystyle (\omega^0 (J\cup\{ 2\},x)\wedge 
\frac{dx_{n-i+1}}{x_{n-i+1}-1} \wedge \omega (T))=
\omega^0 (J,\xi )\wedge \frac{dx_{n-i+1}}{x_{n-i+1}-1} \wedge \omega (T)$.
Therefore $r$ is surjective. Suppose that the element
\begin{equation}
\label{eqn:first kernel}
\eta = \sum_{J\subset [2, n-i]} a_J\omega^0 (J ,x)\wedge 
\frac{dx_{n-i+1}}{x_{n-i+1}-1} \wedge \omega (T)
\end{equation}
is contained in the kernel of $r$.
For a subset $K \subset [2,n-i]$ such that $2 \in K$, the element 
$\displaystyle\tau =
\omega^0 (\partial K ,x)\wedge 
\frac{dx_{n-i+1}}{x_{n-i+1}-1} \wedge \omega (T)$ in
$V_{k,n-i+1, n-i}^{T}$
contains only one term
$\displaystyle\omega^0 (J ,x)\wedge 
\frac{dx_{n-i+1}}{x_{n-i+1}-1} \wedge \omega (T)$ such that $2 \notin J$,
and in this case, $J=K-\{ 2\}$. Therefore there exist an element
$\delta \in V_{k,n-i+1, n-i}^{T}$ such that
$\eta' =\eta - \delta$ contains no term with $2 \notin J$ in expression 
(\ref{eqn:first kernel}).
The relation $r(\eta')=0$ and formula (\ref{eqn:residue})
implies $\eta'=0$.

(2) {\bf The morphism} $r:V_{k,p,n-i+1}^{T} \to V_{k,p,n-i}^{T,\xi}$.
Let  $2 \in K \subset [2,n-i+1]$. 
Then we have
\begin{align}
\label{eqn:second residue of partial}
& r(\omega^0 (\partial K ,x)\wedge \frac{dx_{p}}{x_{p}-1} \wedge \omega (T)) \\
= &
\begin{cases}
\omega^0 (K-\{ 2, n-i+1 \} ,\xi)\wedge \frac{dx_{p}}{x_{p}-1} \wedge \omega (T)
& (\text{ if } n-i+1 \in K) \\
0
& (\text{ if } n-i+1 \notin K), \nonumber
\end{cases}
\end{align}
and the surjectivity follows from this equality. Suppose that the element
\begin{equation}
\label{eqn:second kernel}
\eta = \sum_{J\subset [2, n-i+1]} a_J\omega^0 (\partial J ,x)\wedge 
\frac{dx_{p}}{x_{p}-1} \wedge \omega (T)
\end{equation}
is contained in the kernel of $r$.
Let $2 \in K \subset [2, n-i+1]$. The boundary $\partial K = \sum_J b_J J$
contains only one term with $2 \notin J$.
By using the relation $\sum_{J} b_J\omega^0 (\partial J ,x)\wedge 
\frac{dx_{p}}{x_{p}-1} \wedge \omega (T) =0$, we may assume that 
there is no term with $2 \notin J$ in the expression of (\ref{eqn:second kernel}).
Then by formula (\ref{eqn:second residue of partial}) and $r(\eta)=0$,
there is no term with $n-i+1 \in J$ in the expression of 
(\ref{eqn:second kernel}).
Therefore $\eta \in V_{k,p,n-i}^{T}$.
\end{proof}
We completed the proof of Theorem \ref{thm:de Rham cohom bl up}.1.
The proof of Theorem \ref{thm:de Rham cohom bl up}.2 is similar.
We omit the proof of Theorem \ref{thm:de Rham cohom bl up}.2.

\subsection{Proper transform of $S$-diagonal varieties}

We define divisors $B_i$ in $(\bold A^1)^n$ 
by $B_0 =\{ x_1=0\}$, $B_n=\{ x_n=1\}$ and 
$B_{i}=\{ x_i= x_{i+1}\}$ for $i=1, \dots , n-1$.
Let $S$ be a subset of $[0,n]$ such that $S \neq [0,n]$.
We define $B_S = \cap_{i \in S}B_i$ and $B_{\emptyset}=X_0$.
Then $\dim B_S = \# \bar S -1$, where $\bar S = [0, n]-S$.
Set $\bar S=\{ i_1, \dots ,i_s\}$. Then a coordinate
of $B_S$ is given by $x_{i_2}, \dots x_{i_s}$, since
$$
x_{{i_1}+1}= \cdots = x_{i_2}, 
 \dots, 
x_{{i_{s-1}}+1}= \cdots = x_{i_s}. 
$$

The proper transform of $B_i$ in $Y$ is denoted by $B_i^{pr}$.
Then it is easy to see that 
\begin{align}
\label{first boundary component}
B_0^{pr} \simeq Y(x_2, \dots , x_n).
\end{align}
Let us define divisors $\tilde B_i$ in $Y$ by
\begin{align*}
& \tilde B_0 = B_0^{pr}\cup E_0^{pr}\cup \cdots \cup E_{n-2}^{pr}, \\
& \tilde B_i =B_i^{pr}, \\
& \tilde B_n = B_n^{pr}\cup F_0^{pr}\cup \cdots \cup F_{n-2}^{pr}, 
\end{align*}
where $E_i^{pr}$ (resp. $F_i^{pr}$) is the proper transform of $E_i$
(resp. $F_i$) in $Y$. We put $\tilde B_S = \cap_{i \in S}\tilde B_i$.
We define $B_i^{pr,0} = B_{i}^{pr}\cap Y^0$, 
$\tilde B_i^0= \tilde B_i \cap Y^0$ and 
$\tilde B_S^0 = \tilde B_S \cap Y^0$. 
Then via isomorphism 
(\ref{first boundary component}), we have
$B_0^{pr,0} \simeq Y^0(x_2, \dots , x_n)$.
We can easily check the following lemma.
\begin{lemma}
\label{lemma:exceptioanl and fibration}
\begin{enumerate}
\item
There exists natural homomorphisms
\begin{align}
\label{proper transform blowup}
\begin{cases}
E_i^{pr} \to Y(x_{n-i+1}, \dots ,x_n) \\
F_j^{pr} \to Y(x_1, \dots , x_j),
\end{cases}
\end{align}
which are birational to the morphisms $E_i \to Z_i^{pr}$ and 
$F_j \to W_j^{pr}$ in \S \ref{subsec:successive}.
\item
Let $E_i^{pr,0} = E_i^{pr} \cap Y^0$ (resp. $F_j^{pr,0} = F_j^{pr} \cap Y^0$).
Then the morphisms (\ref{proper transform blowup})
induce morphisms
\begin{align*}
& E_i^{pr,0} \to Y^0(x_{n-i+1}, \dots ,x_n) \\
& F_j^{pr,0} \to Y^0(x_1, \dots , x_j)
\end{align*}
and they are trivial $\bold A^{n-i-1}$-bundle and trivial $\bold A^{n-j-1}$-bundle,
respcetivly. Moreover $\eta_1 = x_1/x_2, \dots , 
\eta_{n-i-1}=x_{n-i-1}/x_{n-i}$ gives a trivialization of 
$\bold A^{n-i-1}$-bundle $E_i^{pr,0} \to Y^0(x_{n-i+1},\dots , x_n)$.
\item
\label{item fiber equation}
The intersections $E_p^{pr}\cap E_i^{pr,0}$ ($p>i$) and 
$B_q^{pr}\cap E_i^{pr,0}$ ( $0 \leq q < n-i-1$)
are horizontal divisors with respect to the morphism
$E_i^{pr,0} \to Y^0(x_{n-i+1}, \dots , x_n)$.
The defining equations of them are 
\begin{align}
\label{fiber equation}
\begin{cases}
E_i^{pr,0}\cap E_{n-2}^{pr}: \eta_2 = 0, \dots ,
E_i^{pr,0}\cap E_{i+1}^{pr}:\eta_{n-i-1}=0, \\
E_i^{pr,0}\cap B_0^{pr}: \eta_1 = 0, 
E_i^{pr,0}\cap B_1^{pr}:\eta_1 = 1, \\
E_i^{pr,0}\cap B_2^{pr}:\eta_2=1,\dots, 
E_i^{pr,0}\cap B^{pr}_{n-i-1}: \eta_{n-i-1}=1.
\end{cases}
\end{align}
\item
\label{pached retraction}
Using the coordinate in (\ref{fiber equation}) of 
Lemma \ref{lemma:exceptioanl and fibration}.\ref{item fiber equation}, the projection 
$$
(\eta_1,\eta_2, \dots, \eta_{n-i-1}) \mapsto (\eta_2, \dots, \eta_{n-i-1})
$$
gives a morphism $\pi_i : E_i^{pr,0} \to E_i^{pr,0} \cap B_0^{pr,0}$.
Moreover the morphisms $\pi_i$ patched together into the map 
$\pi : \tilde B_0^0 \to B_0^{pr,0}$.  The morphism $\pi$ is
homotopy equivalent.
\end{enumerate}
\end{lemma}
\begin{definition}
Via the isomorphism \ref{first boundary component},
the divisor of $Y^0(x_2, \dots ,x_n)$ corresponding to
the divisor $\tilde B_i^0$ is denoted by $\tilde B_i^0(x_2, \dots, x_n)$.
In general, for a set of coordinates $(y_1, \dots ,y_m)$, the
divisor of $Y^0(y_1, \dots , y_m)$ corresponding to $\tilde B_i^0$
is denoted by $\tilde B_i^0(y_1, \dots ,  y_m)$.
Then we have
\begin{align*}
& \tilde B_1^0(x_2, \dots ,x_n) = \cup_{i=0}^{n-2}E_i^{pr} \cap B_0^{pr,0}, \\
& \tilde B_i^0(x_2, \dots ,x_n) = B_i^{pr} \cap B_0^{pr,0}, \\
& \tilde B_n^0(x_2, \dots ,x_n) = \tilde B_n \cap B_0^{pr,0}.
\end{align*}
For a set $S' \subset [1, n]$, we define 
$ \tilde B_{S'}^0(x_2, \dots , x_n)= 
\cap_{i \in S'}\tilde B_i^0(x_2, \dots, x_n)$.
\end{definition}
 The following proposition is a direct consequence of the definition
of $\pi : \tilde B^0_0 \to B^{pr,0}_0$ defined in 
Lemma \ref{lemma:exceptioanl and fibration}.\ref{pached retraction}.
\begin{proposition}
\label{prop:retraction}
Let $S \subset [0,n]$ and $0 \in S$. We put $S'=S-\{ 0\}$.
\begin{enumerate}
\item
$\pi (\tilde B_S^0) = \tilde B_{S'}^0(x_2, \dots , x_n)$.
\item
We have the following commutative diagram, where the horizontal arrows are 
homotopy eqauivalent.

\begin{center}
\vskip 0.1in
\begin{tabular}{ccc}
$\tilde B_0^0$  & $\overset{\pi}\longrightarrow $ & $B_0^{pr,0}$ \\
$\uparrow$ & & $\uparrow$ \\
$\tilde B_S^0$ & $\underset{\text{homotopy equivalent}}\longrightarrow$ &
$\tilde B_{S'}^0(x_2, \dots ,x_n)$  
\end{tabular}
\vskip 0.1in
\end{center}

\end{enumerate}
\end{proposition}

By using Proposition \ref{prop:retraction} successively, we have the
following corollary.
\begin{corollary}
\label{cor:diag homotopy equiv}
Let $S \subset [0,n]$ such that $0 \in S$ and $\bar S = \{ i_1, \dots ,i_s\}$.
The variety $\tilde B_S^0$ is homotopy equivalent to
\begin{align}
\label{sccessive retraction}
E_{n-i_1}^0 \cap B_1^{pr,0} \cap \dots \cap B_{i_1-1}^{pr,0}\cap
& (\cap_{i_1 <j < i_s, j \in S}B_j^{pr,0})    \\
& \cap B_{i_s+1}^{pr,0} \cap \dots \cap B_n^{pr,0}\cap F_{i_s-1}^0
\nonumber
\end{align}
\end{corollary}
Via the isomorphism
$$
E_{n-i_1}^0 \cap B_1^{pr,0} \cap \dots \cap B_{i_1-1}^{pr,0}\cap
B_{i_s+1}^{pr,0} \cap \dots \cap B_n^{pr,0}\cap F_{i_s-1}^0
\simeq Y^0(x_{i_1+1}, \dots ,x_{i_s}),
$$
(\ref{sccessive retraction})
is isomorphic to 
$$
\tilde B_{\hat S}^0(x_{i_1+1}, \dots, x_{i_s})
=\cap_{j\in \hat S}
\tilde B_{j}^0(x_{i_1+1}, \dots, x_{i_s})
\subset Y^0(x_{i_1+1}, \dots , x_{i_s}), 
$$
where 
$\hat S = [i_1, i_s]\cap S$.
$\tilde B_{\hat S}^0(x_{i_1+1}, \dots, x_{i_s})$ is 
identified with $Y^0(x_{i_2}, x_{i_3}, \dots ,x_{i_s})$
and is denoted by $B_S^{pr,0}$.

All the varieties in this paper are defined over $K=\bold Q$.
For a smooth variety $S$ over $K$, the de Rham complex 
$\Omega_{S/K}^{\bullet}$ is denoted by $K_{S,DR}$. 
Then by the definition of de Rham cohomology, we have
$H^i_{DR}(S/K)= H^i(K_{S,DR})$. More generally, for a normal corssing
variety $V=\cup_i V_i$, we define the de Rham complex $K_{V,DR}$
of $V$ by the complex
$$
\oplus_{\# I =1} K_{V_I,DR} \to  \oplus_{\# I =2} K_{V_I,DR} \to
\cdots ,
$$
where $V_I = \cap_{i \in I} V_i$.

\begin{corollary}
The morphism $\tilde B^0_S \to B^{pr,0}_S$ in 
Corollary \ref{cor:diag homotopy equiv} gives the following
quasi-isomorphism
\begin{align}
\label{retract quasi-iso}
K_{B^{pr,0}_S,DR} \to K_{\tilde B^0_S,DR}.
\end{align}
\end{corollary}
Let 
$$
B_S^0 = \{ (x_{i_2}, \dots , x_{i_s}) \mid
x_{i_2} \neq 1, x_{i_3}\neq 0,1, \dots ,x_{i_{s-1}} \neq 0,1,
x_{i_s} \neq 0 \}
$$
Then we have $B_S^0 \subset B_S^{pr,0}$ and
the induced map $H^k(B_S^{pr,0}) \subset
H^k(B_S^0)$ is injective. 
By Theorem \ref{thm:de Rham cohom bl up}, we have the following corollary.

\begin{corollary}
\label{cor:proper trans diag diff}
Under the notation of (\ref{eqn:diff notation}),
the cohomology group $H_{DR}^k(B_S^{pr,0})$ is generated freely by
$$
\omega (j_1, \dots ,j_p;\epsilon_1,\dots, \epsilon_k),
$$
where
$\{j_1 < \dots < j_k\}\subset \{i_2, \dots ,i_s \}$,
$\epsilon_p \in \{0, 1\}$ for $p=2, \dots , k-1$, 
$\epsilon_1 =1$ and $\epsilon_k =0$.
\end{corollary}

\begin{proposition}
\label{commutativity for first differential}
Let $S \subset [0, n]$ and put
$\bar S = \{ i_1, \dots ,i_s\}$. 
We define $S^{(m)} \subset [0,n]$ such that
$\bar S^{(m)} = \{ i_1, \dots ,i_{m-1},i_{m+1},\dots ,i_s\}$.
We put 
\begin{align*}
B_{S^{(0)}}^{00}=& \{ (x_{i_3}, \dots , x_{i_s}) \mid
x_{i_3}\neq 0,1, \dots ,x_{i_{s-1}} \neq 0,1,
x_{i_s} \neq 0 \}, \\
B_{S^{(m)}}^{00}=& \{ (x_{i_2}, \dots , x_{i_{m-1}}, x_{i_{m+1}},
\dots , x_{i_s}) \mid \\
& x_{i_2}\neq 1, x_{i_2} \neq 0,1, \dots ,x_{i_{s-1}} \neq 0,1,
x_{i_s} \neq 0 \}, \\
B_{S^{(s)}}^{00}=& \{ (x_{i_2}, \dots , x_{i_{s-1}}) \mid
x_{i_2}\neq 1,x_{i_3}\neq 0,1, \dots ,x_{i_{s-1}} \neq 0,1 \},
\end{align*}
for $m=1, \dots, s-1$.
Then we have the following commutative diagram for
$m=0, \dots ,s$.

\begin{center}
\vskip 0.1in
\begin{tabular}{ccc}
$B_{S^{(m)}}^{00}$ & $\hookrightarrow$ &
$B_S^0$ 
\\
$\downarrow$ & & $\downarrow$ \\
$B_{S^{(m)}}^{pr,0}$ &  &
$B_S^{pr,0}$ \\
$\uparrow$ & & $\uparrow$ \\
$\tilde B_{S^{(m)}}^{0}$ & $\hookrightarrow$  &
$\tilde B_S^{0}$ 
\end{tabular}
\end{center}
\vskip 0.1in
\end{proposition}

\subsection{Relative de Rham cohomology groups}
\label{subsec:relative de Rham cohom}
We use the same notations $\tilde B_i^0$, $\tilde B_S^0$,  
$B_i^{pr,0}$ and $B_S^{pr,0}$ as in the last subsection.
Then the divisor $\bold B^0 = \cup_{i=0}^n \tilde B_i^0$ of $Y^0$
is normal crossing. The complement of $\bold B^0$ in $Y^0$
is denoted by $V$ and the natural inclusions $V \to Y^0$ and
$\bold B^0 \to Y^0$ are denoted by $j$ and $i$, respectively.
$$
V \overset{j}\to Y^0 \overset{i}\leftarrow  \bold B^0
$$
For subsets $S, T$ of $[0,n]$, such that $S\subset T \neq [0,n]$,
the natural closed immersion $\tilde B_T^0 \to \tilde B_S^0$ is defined.
Let $j_!K_{V,DR}$ be the cone
$Cone(K_{Y^0,DR} \to K_{\bold B^0,DR})$.

Using the natural restriction morphism 
we have a complex
\begin{align}
\label{eqn:complex}
\oplus_{\substack{\# S = * \\ S \subset [0,n]}}  K_{\tilde B_S^0,DR} = 
 (\oplus_{\substack{\# S = 0 \\ S \subset [0,n]}}K_{\tilde B_S^0,DR} & \to   
\oplus_{\substack{\# S = 1 \\ S \subset [0,n]}}K_{\tilde B_S^0,DR}  \\
 \cdots & \to 
\oplus_{\substack{\# S = n \\ S \subset [0,n]}}K_{\tilde B_S^0,DR} ) 
\nonumber
\end{align}
\begin{proposition}
The natural morphism
$$
j_!K_{V,DR} \to \oplus_{\substack{\# S = * \\ S \subset [0,n]}}
K_{\tilde B_S^0,DR} 
$$
is a quasi-isomorphism.
\end{proposition}
We have the spectral sequence 
\begin{equation}
\label{eqn:spect}
E_1^{p,q}=\oplus_{\substack{\# S = p \\ S \subset [0,n]}}
H^{q}(K_{\tilde B_S^{0},DR})
\Rightarrow H^{p+q}(j_!K_{V,DR})
\end{equation}
associated to
the stupid filtration $\sigma^*$:
\begin{align*}
\sigma^i(\oplus_{\substack{\# S = * \\ S \subset [0,n]}}
K_{\tilde B_S^{0},DR}) = 
(0 & \to \cdots \to 0 \to \oplus_{\substack{\# S = i \\ S \subset [0,n]}}
K_{\tilde B_S^{0},DR} \\
\cdots & \to \oplus_{\substack{\# S = n \\ S 
\subset [0,n]}}K_{\tilde B_S^{0},DR}\to 0 )
\end{align*}
By quasi-isomorphism (\ref{retract quasi-iso}),
$E_1$-term of spectral sequence (\ref{eqn:spect}) is
isomorphic to
$E_1^{p,q}=\oplus_{\substack{\# S = p \\ S \subset [0,n]}}
H^{q}_{DR}(B_S^{pr,0})$. To compute the differential $d_1$ of
the spectral sequence, we use the following commutative
diagram, which is a consequence of 
Proposition \ref{commutativity for first differential}.

\begin{center}
\vskip 0.1in
\begin{tabular}{ccc}
$H^q_{DR}(B_S^{pr,0})$ & $\rightarrow$ &
$H^q_{DR}(B_{S^{(m)}}^{pr,0})$ \\
$\downarrow$ & & $\downarrow$ \\
$H^q_{DR}(B_S^{0})$ & $\rightarrow$ &
$H^q_{DR}(B_{S^{(m)}}^{0,0})$, 
\end{tabular}
\end{center}
\vskip 0.1in
where $S, S^{(m)}$ are as in 
Proposition \ref{commutativity for first differential}.

We define the word type of 
$\omega (j_1, \dots ,j_p;\epsilon_1,\dots, \epsilon_k) $
in (\ref{eqn:diff notation})
as $(\epsilon_1,\dots, \epsilon_k)$. 
For a word $W=(\epsilon_1,\dots, \epsilon_k)$, $k$
is called the length of $W$, which we denote by $\len (W)$.
Let $\bold W_q$ be the set of words $W=(\epsilon_1, \dots ,\epsilon_q)$
of length $q$ such that $\epsilon_i=0,1$, $\epsilon_1 = 1$
and $\epsilon_q =0$.
Note that $\bold W_0$ consists of the empty word $()$.
\begin{theorem}
\label{thm:purity}
For spectral sequence (\ref{eqn:spect}), we have
$E_2^{p,k}=0$ if $p+k\neq n$ and
$$
E_2^{n-k,k}\simeq \oplus_{W\in \bold W_k}K.
$$
As a consequence, we have $H^i(j_!K_{V,DR})=0$ for $i \neq n$ and
the filtration on $H^n(j_!K_{V,DR})$
induced by spectral sequence (\ref{eqn:spect}) coincides
with the weight filtration.
\end{theorem}
\begin{proof}
The subspace of $H^k_{DR}(B_S^{pr,0})$ generated by differential
forms of word type $W=(\epsilon_1,\dots, \epsilon_k)$
is denoted by $H_{S,W}$. Then, by 
Corollary \ref{cor:proper trans diag diff}, we have 
$$
H^k_{DR}(B_S^{pr,0}) \simeq \oplus_{W \in \bold W_k} H_{S,W}.
$$
We can choose the local coordinate of $B_S^0$ as $(x_{i_2}, \dots , x_{i_s})$
if $\bar S = \{i_1< \cdots < i_s\}$. 
Then the differential $d_1$ of the spectral sequence $E_1^{p,k}$ of 
(\ref{eqn:spect}) preserves the word type and as a consequence, the complex
$E_1^{*,k}$ is a direct sum $\oplus_{W \in \bold W_k}H_W^*$
of the complices $H_W^*$, where
$$
H_W^*: H_W^0 \simeq 
\oplus_{\substack{\# \bar S = n+1 \\ S \subset [0,n]}}H_{S,W} \to 
H_W^1 \simeq \oplus_{\substack{\# \bar S = n \\ S \subset [0,n]}}
H_{S,W} \to \cdots
$$
Here we put $\bar S = [0,n] - S$.
Therefore it is enough to prove the following proposition.
\end{proof}
\begin{proposition}
Let $W \in \bold W_k$.
We have $H^i(H_W^*)=0$ for $i \neq n-k$ and $H^{n-k}(H_W^*) \simeq K$.
\end{proposition}
\begin{proof}
Let $\bar T,\bar S$ be subsets of $[0,n]$ 
such that $\# \bar T=s-1, \# \bar S = s$ and $\bar T \subset \bar S$.
Set $\bar S = \{i_1, \dots, i_s\}$ and $\bar T= \bar S-\{i_l\}$. 
We put $S=[0,n]-\bar S$ and $T=[0,n]-\bar T$.
Then the restriction 
$i_{TS}:B_T^{pr,0} \to B_S^{pr,0}$ is given as follows.
\begin{enumerate}
\item
If $l\neq 1, s$, 
$i_{TS}^*(x_{i_p})= x_{i_p}$ if $p\neq l$ and 
$i_{TS}^*(x_{i_l})=x_{i_{l+1}}$.
\item
If $l = 1$,
$i_{TS}^*(x_{i_p})=x_{i_p}$ if $p\neq 2$ and
$i_{TS}^*(x_{i_2})=0$.
\item
If $l = s$,
$i_{TS}^*(x_{i_p})=x_{i_p}$ if $p\neq s$ and
$i_{TS}^*(x_{i_s})=1$.
\end{enumerate}
For an integer $0 \leq i \leq n+1$
and $W \in \bold W_k$, we define the filtration $F^i(H_W^p)$ of
$H_W^p$ to be the subspace generated by
$H_{S,W}$ for $\bar S \subset [i,n]$ and $\# \bar S = n+1-p$.
Then we can check that the restriction $i_{TS}\*$ preserves the
filtration and $F^*$ is a filtration of complex $H_W^*$.
The differential on $Gr_F^i(H_W^*)$ induced by $d_1$
is also denoted by $d_1$. It is enough to prove the following
claim for the proof of the proposition.

{\bf Claim}
Let $W \in \bold W_k$. (1)If $H^j(Gr_F^i(H_W^*))\neq 0$ then $i=j=n-k $,
and (2)$H^{n-k}(Gr_F^{n-k}(H_W^*))=K$.

We put $s=n+1-p$.
We identify $Gr_F^i(H_W^{p})$ with the space 
$\oplus_{\bar S \in\Cal S_{i,s}} H_{S,W}$,
where 
$$
\Cal S_{i,s} = \{\bar S=\{i<i_2<\cdots < i_s\}\mid \# \bar S= s\}.
$$
The set $\Cal S_{i,s}$ is identified with 
$\{\bar S'\mid \bar S'\subset [i+1,n], \# \bar S'= s-1\}$.
The class of
$\omega (j_1, \dots ,j_k,W) \in H_{S,W}$ in $Gr_F^i(H_W^{p})$ is denoted by
$(\bar S', \tau)$, where 
$\bar S'=\{i_2, \dots ,i_s\}$ and $\tau=(j_1, \dots ,j_k)$. 
The set 
$$
\{ (\bar S',\tau ) \mid
\bar S'\subset [i+1,n],\tau \subset \bar S',
\# \bar S' = s-1 = n-p, \#\tau = k \}
$$
forms a base of $Gr_F^i(H_W^p)$.
Let $\bar T'$ be an element of $\Cal S_{i,s-1}$ such that 
$\bar T' \subset \bar S'$.
The element $\bar S' -\bar T'$ is denoted by $i_l$.
Then the $\bar T'$ component of $d_1(\bar S', \tau)$ is 
\begin{enumerate}
\item
$(\bar T',\tau)$ if $i_l \notin \tau$, 
\item
$(\bar T',\tau')$, where $\tau' = 
(j_1, \dots ,j_{m-1}, i_{l+1}, j_{m+1}, \dots , j_k)$
if $i_l=j_m $, $l\neq s$ and $i_{l+1} \notin \tau$,
\item
zero if $l \neq s$,  $i_l \in \tau$, and $i_{l+1} \in \tau$, and
\item
zero if $l=s$.
\end{enumerate}
We introduce a partial order on the set 
$\Cal T =\{ \tau \subset [i+1,n]\mid \# \tau =k\}$
so that
$(j_{11}, \dots ,j_{1k}) \leq (j_{21}, \dots ,j_{2k})$ if and only if 
$j_{1m} \leq j_{2m}$ for $m=1, \dots ,k$. We choose a numbering 
$\Cal T =  \{\tau_1, \dots , \tau_{\kappa}\}$ on $\Cal T$
such that $\tau_i \leq \tau_j$ implies $i \leq j$.
We define $G^t(Gr_F^i(H_W^p))$ to be the subspace of $Gr_F^i(H_W^p)$ 
generated by $(\bar S', \tau_j)$ for $t \leq j$. Then $G^*$ defines
a decreasing filtration on $Gr_F^i(H_W^p)$ preserved by the 
differential $d_1$.  It is easy to see that the complex $Gr_G^tGr_F^i(H_W^*)$
is exact if $n-i > k$. If $n-i=k$, then $Gr_F^i(H_W^p)=0$ if $p \neq n-i$ and
$Gr_F^i(H_W^{n-i})$ is a one dimensional vector space 
generated by $(\bar S' , \tau) =([i+1,n],[i+1,n])$.
Therefore we have the claim.
\end{proof}

\section{Generators of mixed Tate motives}
\label{sec:generator MTM}
\subsection{Review of results on mixed Tate motives}

In this subsection, we recall the theory of mixed Tate motives
and their properties.
M.Levine constructed the derived category $DTM_{\bold Q}$
of mixed Tate motives in \cite{L}.
This is an additive category with the following properties.

(1) {\bf Triangulated category.} (See \cite{L} p.19 Definition 2.1.6., 
p.45 Corollary 3.4.3.) $DTM_\bold Q$ is a
triangulated category. The set 
of distinguished triangles and the shift operator $A \mapsto A[1]$ are specified.
A shift operator is an equivalence of category.
A triangle is a pair of morphisms and objects
$A \to B \to C \to A[1]$, which is denoted by $A \to B \to C \overset{+1}\to$
for short. The following axioms are imposed.
\begin{enumerate}
\item
For a morphism $f: A \to B$, there exists a unique object and morphisms
$g$ and $\delta$ such that
$C \overset{g}\to A \overset{f} \to B \underset{+1}{\overset{\delta}\to}$
is a distinguished triangle.
\item
If $A\to B\to C \overset{+1}\to$ is a distinguished triangle, then 
$B\to C \to A[1]\overset{+1}\to$ and $C[-1]\to A \to B \overset{+1}\to$
are distinguished triangles.
\item
Let $A \to B \to C \overset{+1}\to$ be a distinguished triangle and $X$ be an object of 
$DTM_\bold Q$. Then we have the following long exact sequences:
$$
\cdots \to Hom(X, A) \to Hom(X, B) \to Hom(X, C) \to Hom(X, A[1]) \to \cdots
$$
$$
\cdots \to Hom(C, X) \to Hom(B, X) \to Hom(A, X) \to Hom(C, X[1]) \to \cdots
$$
\end{enumerate}

(2) {\bf Geometric objects.} 
\begin{enumerate}
\item
For an object $A$ in $DTM_\bold Q$ and an integer $k$, 
the Tate twist $A(k)$ is defined.
\item
An object $\bold Q_{X}$ of $DTM_\bold Q$ is attached
to a variety $X$ over $\bold Q$ with a stratification
$X_0 = X \supset X_1 \supset X_2 \cdots$ satisfying
$$
X_i - X_{i-1} = \coprod_{\text{finite}} \bold A^{m-i},
$$
Moreover it is
contravariant with respect to $X$.
\item
Let $X_1$ and $X_2$ be smooth varieties as above and $X_1 \to X_2$ be
a closed embedding of codimension $d$. Set $U = X_2 - X_1$. Then $\bold Q_U$ in
$DTM_\bold Q$ is defined functorially on $U$. Moreover there
exists a morphism $\bold Q_{X_1}(d)[2d] \to \bold Q_{X_2}$
such that 
$\bold Q_{X_1}(d)[2d] \to \bold Q_{X_2} \to \bold Q_U \overset{+1}\to$
is a distinguished triangle.
\end{enumerate}

(3) {\bf Hodge realization.}(See \cite{L} \S 2.3, p.273.)
There exists a realization functor 
$\bold H :DTM_\bold Q \to GMHS$ from $DTM_\bold Q$
to the category of graded objects of mixed Hodge structures with the 
following properties.
The degree $i$-part of $\bold H$ is denoted by $\bold H^i$.
\begin{enumerate}
\item
For a variety $U$ as in (2), $\bold H^i(\bold Q_U) = 
(H^i_B(U, \bold Q),H^i_{DR}(U, \bold Q))$.
\item
For a distinguished triangle $A \to B \to C \overset{+1}\to$, we have 
a long exact sequence:
$$
\cdots \to \bold H^i(A) \to \bold H^i(B) \to 
\bold H^i (C) \to \bold H^{i+1} (A) \to \cdots
$$
\end{enumerate}

(5){\bf Relations to the K-groups.}
The following equality holds:
\begin{align*}
Hom_{DTM_\bold Q}(\bold Q, \bold Q(i)[1]) & = K_{2i-1}(\bold Q)
\otimes \bold Q, 
\text{ for }i=1, 2, \dots \\
Hom_{DTM_\bold Q}(\bold Q, \bold Q(i)[2]) & = 0. 
\end{align*}

(4){\bf Weight filtration.}
We consider objects with weight filtration.
Weight filtration of a object $A$ is a sequence of morphisms
$W_iA \to A$, $W_iA \to W_{i+1}A$ with the following properties. 
\begin{enumerate} 
\item
For a sufficiently big $i$ (resp. small $i$), $W_iA \to A$
is an isomorphism (resp. zero map).
\item
The cone of $W_{i-1}A \to W_iA$ is an object 
of $DTM_\bold Q$ generated
by $\bold Q(-i)$ in $DTM_\bold Q$.
\end{enumerate}
\begin{remark}
In this paper, $W_i$ denotes the $2i$-th weight filtration
to simplify the notation.
\end{remark}

An object $A$ is called an abelian object if $Cone(W_{i-1} \to W_i)[1]$
is isomorphic to $\bold Q(-i)^{r_i}$ for all $i$.
According to \cite{L2}, the full subcategory $A_{TM}$ of abelian object in 
$DTM_\bold Q$ is an abelian category.

By the definition of $W_i$ and the compatibility with the realization
functor, we have $W_i\bold H^i(A) = \bold H(W_iA)$.
Therefore $A$ is an abelian object if and only if $\bold H^i(A) = 0$
if $i \neq 0$.
The following proposition is a direct consequence of 
Theorem \ref{thm:purity}.
\begin{proposition}
Let $j_!\bold Q_V$ be the cone 
$Cone (\bold Q_{Y^0} \to \bold Q_{\bold B^0})$. Then $j_!\bold Q_V[n]$
is an abelian object.
\end{proposition}

\subsection{Comparison of the extension groups}
\label{sebsec:comparison extension grp}
In this subsection, we compare extension group of
mixed motives and that of mixed Hodge structures for certain
Tate structures.
We define a group homomorphism 
$$
ch : Ext^1_{MHS}(\bold Z, \bold Z(1)) \overset{\simeq}\to \bold C^{\times}
$$
as follows. Let
\begin{equation}
\label{eqn:ext class}
u_{MHS}:0\to \bold Z(1) \to M_{\bold Z} \overset{\pi}\to \bold Z \to 0
\end{equation}
be an extension of mixed Hodge structures. Set 
$M_{\bold C} =M_{\bold Z} \otimes \bold C$.
Then the natural map $F^0(M_{\bold C}) \to \bold Z\otimes \bold C$
is an isomorphism. Let $\omega$ be the element in $F^0(M_{\bold C})$
corresponding to $1$ via this isomorphism. 
We consider the dual of exact sequence (\ref{eqn:ext class}):
$$
u_{MHS}^*:0\to \bold Z \to M_{\bold Z}^* \overset{\pi'}\to \bold Z(-1) \to 0
$$
Let $\gamma$ be the element in $M_{\bold Z}^*$ such that 
$\pi'(\gamma)=2\pi i$. 
Put $ch (u_{MHS})=\exp (\langle\omega, \gamma\rangle)\in \bold C^{\times}$. 
Then it is easy to see that $ch (u_{MHS})$
does not depend on the choice of $\gamma$.

Let $CH^1(\bold Q, 1)$ be the Bloch's higher Chow group. (See \cite{B}.)
In \cite{L}, Levine defined an isomorphism
$cl : CH^1(\bold Q,1)\otimes\bold Q \to 
Hom_{DTM_\bold Q}(\bold Q, \bold Q(1)[1])$,
which is called the cycle map. By using the isomorphism 
\begin{equation}
\label{eqn:higher chow}
\bold Q^{\times} \simeq CH^1(\bold Q,1),
\end{equation}
we have an isomorphism
$\bold Q^{\times}\otimes \bold Q \to Hom_{DTM_\bold Q}(\bold Q, \bold Q(1)[1])$,
which is also denoted by $cl$.
We have the following proposition
\begin{proposition}
\label{prop:compatiblity ext one}
The following diagram commutes.

\begin{tabular}{cccl}
$Hom_{DTM_\bold Q}(\bold Q, \bold Q(1)[1])$ & 
$\underset{\simeq}{\overset{cl}\leftarrow}$ & 
$\bold Q^{\times}\otimes \bold Q$ & \\
$\bold H \downarrow$&   & $\downarrow$ & natural inclusion   \\
$Ext^1_{MHS}(\bold Q,  \bold Q(1) )$& $\underset{\simeq}\to$ & 
$\bold C^{\times}\otimes\bold Q$ &\\
\end{tabular}
\end{proposition}
We recall the definition of the cycle map. (See \cite{L}.)
Let $\Delta^p$ be the variety over $\bold Q$ defined by
$$
\Delta^p =\{ (x_0, \dots ,x_p) \mid \sum_{i=0}^p x_i =1 \}.
$$
For a subset $S$ of $[0,p]$, the subvariety of $\Delta^p$
defined by $x_i = 0$ for $i\in S$ is denoted by $\delta_S$.
Let $z^q(\bold Q, p)$ be the free $\bold Z$ module generated by
codimension $q$ cycles $\gamma$ of $\Delta^p$, where the codimension
of $\gamma\cap\delta_S$ is at least $q$.
Since $\{ \Delta^p \}$ is a cosimplicial scheme, we have a 
complex $z^q(\bold Q,*)$ using restrictions to faces.
The homology of $H_p(z^q(\bold Q, *))$ is denoted by $CH^q(p)$
and it is called the Bloch's higher Chow group.
An element in $CH^1(1)\otimes \bold Q$ is represented 
by a $\bold Q$-linear combination of
0-dimensional sub-schemes $z_1, \dots ,z_j$ 
in $\Delta^1$ defined over $\bold Q$ which does not intersect
with $\delta_{\{0\}}\cup \delta_{\{1\}}$. 
Equality (\ref{eqn:higher chow})
is obtained by attaching the class of $z_i$ to
$Nm(z_i-1)/Nm(z_i) \in \bold Q^{\times}$. 
By using trace, $z_i$ defines a 
morphism $[z_i]:\bold Q \to \bold Q_{\Delta^1}(1)[2]$.
After \cite{L}, we define $\bold Q_{\Delta^*}^{\leq 1}(1)[2]$
by
$$
\bold Q_{\Delta^*}^{\leq 1}(1)[2]=
(\bold Q_{\Delta^1}(1)[2]_{\text{degree}=-1} \to 
\bold Q_{\Delta^0}(1)[2]^{\oplus 2}_{\text{degree}=0}).
$$
Then $\bold Q_{\Delta^*}^{\leq 1}(1)[2]$ is quasi-isomorphic to $\bold Q(1)[2]$.
The morphism 
$[z_i]:$ 
\linebreak
$\bold Q_{\text{degree}=-1} \to 
\bold Q_{\Delta^1}(1)[2]_{\text{degree}=-1}$
defines a morphism $\bold Q[1] \to \bold Q_{\Delta^*}^{\leq 1}(1)[2]$
and this defines a morphism $cl(z_i) \in Hom_{DTM}(\bold Q, \bold Q(1)[1])$.
We compute $ch(\bold H(cl(z_i)))$, where $\bold H (cl(z_i))$ is the Hodge
realization of $cl(z_i)$.
For simplicity, we assume $z \in \Delta^1-\{ 0,1\}$ is
a $\bold Q$ rational section. Let $j$ be the inclusion
$\Delta^1-\{ 0,1\}\to \Delta^1$.
Then $z$ induces a morphism of mixed Hodge complex:
$$
[z]:\bold Q \to 
\bold R\Gamma (\Delta^1,j_!\bold Q(1)_{\Delta^1-\{ 0,1\}})[2]\simeq \bold Q(1)[1].
$$
Then we have the triangle
\begin{align*}
\bold Q(1)[1]\simeq \bold R\Gamma (\Delta^1,j_!\bold Q(1)_{\Delta^1-\{ 0,1\}}[2])
& \to
\bold R\Gamma (\Delta^1-\{ z\},j_!\bold Q(1)_{\Delta^1-\{ 0,1\}})[2] \\
& \to
\bold Q[1]
\end{align*}
By taking $H^{-1}$ of complices, we have the exact sequence of
mixed Hodge structures:
\begin{equation}
\label{eqn:Hodge real of cycle}
0 \to
H^1 (\Delta^1,j_!\bold Q(1)_{\Delta^1-\{ 0,1\}})
\to
H^1 (\Delta^1-\{ z\},j_!\bold Q(1)_{\Delta^1-\{ 0,1\}})
\overset{\pi}\to \bold Q \to 0.
\end{equation}
This exact sequence corresponds to the Hodge realization of $cl(z)$.
An element $\frac{1}{2\pi i}\frac{dx}{x-z}$ in
$H^1 _{DR}(\Delta^1-\{ z\},j_!\bold Q(1)_{\Delta^1-\{ 0,1\}})$
corresponds to $1 \in \bold Q$ via the projection $\pi$ in 
(\ref{eqn:Hodge real of cycle}).
The cycle $2\pi i\cdot [0,1]$ defines an element
$H_1(\Delta^1-\{ z\}, \text{ mod }\{ 0,1\})(1)$,
whose image in
$H_1(\Delta^1, \text{ mod }\{ 0,1\})(1)$ is $2 \pi i$.
Then we have
\begin{align*}
ch(\bold H (cl(z))) & = \exp(2\pi i\int_0^1\frac{1}{2\pi i}\frac{dx}{x-z}) \\
& = \frac{z-1}{z}
\end{align*}
Therefore the diagram in the proposition commutes.
\subsection{Splitting in level $(a_i)$}

Let $A$ be an abelian object in $\Cal A_{TM}$ and $S=\{ a_1,\dots ,a_k \}$ be 
a set of integers such that $a_i  < a_{i+1}$.
An abelian object $A$ in $A_{TM}$ is said to be of type $S$ if
$Gr_a^WA = 0$ for $a \notin S$ and the full subcategory of abelian objects
in $\Cal A_{TM}$ of type $S$ is denoted by $\Cal A_S$.
Note that the category $\Cal A_S$ is stable under taking 
direct sums and subquotients.
Let $A$ be an object of type $S$ and $a=a_i, b=a_{i+1} \in S$. We have 
the following exact sequence in $\Cal A_{TM}$:
$$
0 \to W_aA/W_{a-1}A \to W_bA/W_{a-1}A \to W_bA/W_{b-1}A \to 0.
$$
If the morphism $\phi :Hom_{DTM_\bold Q}(W_{b}A/W_{b-1}A, W_aA/W_{a-1}A[1])$
corresponding to the above exact sequence is zero,
$A$ is said to split in level $(a,b)$.
If $A$ splits in level $(a,b)$ then a subquotient of $A$ splits 
in level $(a, b)$.
It is easy to see that a objects of type $S-\{a_i\}$ and $S-\{a_{i+1}\}$
split in level $(a, b)$.
\begin{proposition}
\label{prop: decomp}
Let $A$ be an object of type $S$ which splits
in level $(a_i, a_{i+1})$.
Then $A$ is a subquotient of a module $B = B_1\oplus B_2$,
where (1) $B_1$ (resp. $B_2$) is of type $S-\{a_{i+1}\}$
(resp. $S-\{a_{i}\}$), 
(2)$W_{a_{i}}B_1$ (resp. $W_{a_{i-1}}B_2$) 
is isomorphic to the direct sum of copies of
$W_{a_{i}}A$ (resp. $W_{a_{i-1}}A$) and
(3)$B_1/W_{a_{i+1}}B_1$ (resp. $B_1/W_{a_{i}}B_2$) 
is isomorphic to the direct sum of copies of
$A/W_{a_{i+1}}A$ (resp. $A/W_{a_{i}}A$).
Here we use the notation $W_{a_0}A=0$ and $W_{a_{k+1}}A=A$.
\end{proposition}

Since $Hom_{DTM_\bold Q}(X, Y[2])=0$ for abelian objects $X, Y$,
the morphism 
$$
Hom_{DTM_\bold Q}(X_1, Y[1]) \to 
Hom_{DTM_\bold Q}(X_2, Y[1])
$$ 
induced by an injective
morphism $X_2 \to X_1$ is surjective.

We put $a=a_i, b=a_{i+1}$.
By the assumption of the proposition, the exact sequence
$$
0 \to W_aA/W_{a-1}A \to W_{b}A/W_{a-1}A \to W_bA/W_{b-1}A \to 0
$$
splits and we have $W_{b}A/W_{a-1}A =W_aA/W_{a-1}A \oplus
W_bA/W_{b-1}A$.
By pushing forward the exact sequence
$$
0 \to W_bA/W_{a-1}A \to A/W_{a-1}A \to A/W_{b} \to 0
$$
by the morphism $W_{b}A/W_{a-1}A \to W_aA/W_{a-1}WA$, we have the
following exact sequence.
$$
0 \to W_aA/W_{a-1}A \to C_1 \to A/W_{b} \to 0
$$
It is easy to see that the morphism $i: A/W_{a-1} \to C_1 \oplus A/W_{b-1}$
is an injective morphism. Therefore one can find an extension
\begin{equation}
\label{eqn:1}
0 \to W_{a-1}A \to C_2 \to C_1 \oplus A/W_{b-1}\to 0
\end{equation}
such that the pull back by the morphism $i$ is isomorphic to
$$
0 \to W_{a-1}A \to A \to A /W_{a-1} \to 0.
$$
Then $A \to C_2$ is injective.
By pulling back exact sequence (\ref{eqn:1}) by morphisms
$C_1 \to C_1 \oplus A/W_{b-1}$ and
$A/W_{b-1} \to C_1 \oplus A/W_{b-1}$, we have the following exact sequences.
$$
0 \to W_{a-1}A \to B_1 \to C_1 \to 0,
$$
$$
0 \to W_{a-1}A \to B_2 \to A/W_{b-1}\to 0.
$$
Then the morphism $B_1 \oplus B_2 \to C_2$ is surjective.
By pulling back exact sequence (\ref{eqn:1})
by the composite morphism $W_aA/W_{a-1}A \to A/W_{a-1} \to
C_1 \oplus A/W_{b-1}$, We have an exact sequence
$$
0 \to W_aA \to B_1 \to A/W_bA \to 0.
$$
Thus we have the proposition.

\subsection{Generators of $\Cal A_{TM}$}
\label{subsec:generator}

In this section we denote 
$Hom_{\Cal A_{TM}}(A, B[1])$ by
$Ext^1(A, B)$.
We define a subset $\Cal N_n$ by
\begin{align*}
\Cal N_n = \{ \{a_1, \dots , a_n\} \mid
&\ 
0 \leq a_1 < \cdots < a_n \leq n
\text{ and } a_{i+1}-a_i \text{ is odd} \\
& \text{and greater than 1} \}.
\end{align*}
For $S \in \Cal N_n$,
we inductively define a series of objects $\{ M_S \}_{S \subset \Cal N}$,
where $M_S \in \Cal A_S$
and morphisms $M_S \to \bold Q(a_n)$
as follows.
We set $T = \{ a_1, \dots , a_{n-1}\}$. 
\begin{enumerate}
\item
Let $M_{\{ a_1\}}= \bold Q(a_1)$.  
\item
$M_S$ is defined by the extension
$$
0 \to M_T \to M_S \to \bold Q (a_n) \to 0
$$
corresponding to $u \in Ext^1(\bold Q(a_n), M_T)$ such that
$\pi (u) \neq 0$ where $\pi$ is the natural map
$$
Ext^1(\bold Q(a_n), M_T) \to
Ext^1(\bold Q(a_n), \bold Q(a_{n-1}))
$$
induced by $M_T \to \bold Q(a_{n-1})$.
\end{enumerate}
By the construction, $Gr_{a_i}M_S \simeq \bold Q(a_i)$ for $i \in S$.
Let $M \in \Cal A_S$. We define $<M, \Cal A_{T_i}>_{T_i \subsetne S}$
as the minimal sucategroy of $\Cal A_S$ containing $\Cal A_{T_i}$ 
($T_i \subsetne S$) and
$M$ and stable under taking direct sums and subquotients.
\begin{proposition}
\label{prop:generated by one obj}
For $S \in \Cal N_n$, 
we have $\Cal A_S = < M_S, \Cal A_{T_i}>_{T_i \subsetne S}$.
\end{proposition}
We omit the proof of the next lemma.
\begin{lemma}
Let $u_1 , u_2 \in Ext^1(A, B)$ corresponding to the following extensions:
\begin{align*}
0 \to B \to M_1 \to A \to 0, \\
0 \to B \to M_2 \to A \to 0.
\end{align*}
Let $M_3$ be defined by the extension
$$
0 \to B \to M_3 \to A \to 0 
$$
corresponding to the element $u_1 + u_2 \in Ext^1(A, B)$. Then 
$M_3$ is a subquotient of $M_1 \oplus M_2$.
\end{lemma}
\begin{proof}[Proof of Proposition \ref{prop:generated by one obj}]
We prove the proposition by the induction on $S$.
Let $M$ be an element of $\Cal A_S$. We may assume that
$Gr_{a_n}M \simeq \bold Q(a_n)$.
Let $S = \{a_1, \dots , a_n\}$ ($a_1 < \cdots < a_n$) 
and $T = \{ a_1, \dots , a_{n-1}\}$.
Then $L = W_{a_{n}-1} M \in \Cal A_T$. By the hypothesis of the
induction, $L$ can be obtained by the subquotient of 
$M_T^{\oplus m} \oplus\oplus_i L_i$, where $L_i \in \Cal A_{U_i}$ and
$U_i \subsetne T$.
That is, there exists an object $N \in \Cal A_T$ and surjective and
injective morphisms $M_T^{\oplus m} \oplus\oplus_i L_i \to N$
and $L \to N$.
Since $Ext^1(\bold Q(a_n),M_T^{\oplus m}\oplus \oplus_i L_i)
\to Ext^1(\bold Q(a_n), N)$ is surjective,
the object $M$ is a subquotient of the object
$M_1$ defined by 
the exact sequence
$$
0 \to M_T^{\oplus m} \oplus\oplus_i L_i \to M_1
\to  \bold Q(a_n) \to 0.
$$
Suppose that $M_1$ is an object corresponding to
$v\in Ext^1(\bold Q(a_n),M_T^{\oplus m} \oplus\oplus_i L_i)$.
Let us write $v=(v_{11}, \dots ,v_{1m}, v_{21}, \dots ,v_{2k})$,
where $v_{1i} \in $ 
\linebreak
$Ext^1(\bold Q(a_n), M_T)$ and
$v_{2j} \in Ext^1(\bold Q(a_n), L_j)$. 
Let $u \in Ext^1(\bold Q (a_n), M_T)$ be the 
element corresponding to $M_S$.

(1)
Since $L_i \in \Cal A_{U_i}$ and $U_i \subsetne T$,
the extension $A_i$ corresponding to 
$v_{2j} \in Ext^1(\bold Q(a_n), L_j)$ is an element in
$\Cal A_{U_i\cup \{a_n \}}$. Note that 
$U_i\cup \{a_n \}\subsetne S$. 

(2)
By the definition of $M_S$, the image $\pi (u)$ of $u$
is a generator 
of 
\linebreak
$Ext^1(\bold Q (a_n), \bold Q (a_{n-1}))$.
Thus there exist rational numbers $k_1, \dots, k_m$ such
that $\pi (v_{1i}) = k_i \pi (u)$.
Therefore an extension $B_i$ corresponding to 
$v_{1i}-k_i u \in Ext^1(\bold Q (a_n), M_T)$ 
splits in level $(a_{n-1}, a_n)$.
Therefore $B_i$ is a object in
$<\Cal A_{T_i}>_{T_i \subsetne S}$ by Proposition \ref{prop: decomp}.

(3) The extension $C_i$ corresponding to 
$k_i u \in Ext^1(\bold Q(a_n),M_T)$ is isomorphic to 
(1)
$M_S$ if $k_i \neq 0$, and
(2)
$M_T \oplus \bold Q(a_n)$ if $k_i = 0$.

Since $M_1$ is a subquotient of $A_i, B_i, C_i$ as in (1),(2) and (3),
we have the proposition by induction.
\end{proof}

\section{Periods of mixed Hodge structures and multiple zeta values}
\label{sec:periods MTS MZV}
\subsection{Subspace generated by periods}
Let $n \geq 0$ be a natural number and 
$H=(H_B, H_{DR})$ a $\bold Q$-Hodge structre over $\bold Q$
such that $W_{-1}H=0$ and $W_nH = H$.
(Note that $W_i$ denote the $2i$-th weight 
filtration in the usual convention.)
Let $H^*=(H_B^*,H_{DR}^*)$ be the dual of $H$.
We define the period space $p_n(H)$ of $H$
of weight $n$ by the $\bold Q$-linear hull of
the set $\{\langle\gamma, \omega\rangle \mid 
\gamma \in H_B^*, \omega \in F^nH_{DR}\}$.
Then $p_n(H)$ is a finite dimensional $\bold Q$-subvector space in $\bold C$.
The following properties are straight forward from the definition.
Let $H_1,H_2$ be mixed Tate Hodge strucutres 
such that $W_{-1}H=0$ and $W_nH = H$.
\begin{enumerate}
\item
If $H_1$ is a subquotient of $H_2$, we have
$p_n(H_1) \subset p_n(H_2)$.
\item
If $W_{n-1}H_1 = H_1$, then $p_n(H_1) = 0$.
\item
We have $p_{n+k}(H(-k)) =(2\pi i)^{k} p_n(H)$.
\end{enumerate}

For an abelian object $M\in \Cal A_{TM}$, the period space $p_n(\bold H(M))$
of the Hodge realization of $M$ is also denoted by $p_n(M)$.
Since the differnetial form $\omega \in H_{DR}$ is defined over $\bold R$,
the complex conjugation $c$ for tolopolgical realizations and the
complex conjugation 
for the periods $\langle\gamma , \omega \rangle$ commultes, i.e.
$\overline{\langle\gamma , \omega \rangle}=\langle\gamma^c, \omega \rangle$.
As a consequence, $p(H)$ is stable under the action of complex conjugation.
Moreover, the topological action of complex conjugate $c$ acts on
the Betti realization $\bold Q (k)$ as $(-1)^k$-multiplication. 
\begin{lemma}
Let $H$ be an object of $\Cal A_S$, where $S=\{ a_1< \cdots < a_k \}$.
Then the spaces $p_{a_k}(H)$ and $p_{a_k}(H/W_{a_1})$ are stable under
the complex action. The action of complex conjugate on the space 
$p_{a_k}(H)/p_{a_k}(H/W_{a_1})$ coincides with the $(-1)^{a_1}$-multiplicaltion. 
\end{lemma}
\begin{proof}
The space $p_{a_k}(H)$ (resp. $p_{a_k}(H/W_{a_1})$) is generated by the set
$\{ \langle\gamma , \omega\rangle \mid \omega \in F^{a_k}H_{DR},\gamma \in H^*\}$
(resp. $\{ \langle\gamma , \omega\rangle 
\mid \omega \in F^{a_k}H_{DR},\gamma \in W_{-a_2}H^*\}$).
For an element $\gamma \in H_B^*$, we have 
$\gamma -(-1)^{a_1} \gamma^c \in W_{-a_2}H_B^*$. Therefore
$\langle\gamma , \omega\rangle- (-1)^{a_1}
\overline{\langle\gamma , \omega\rangle} \in p(H/W_{a_1})$.
\end{proof}

Let $S$ be a finite subset of $[0,n]$. Let $p_n(S)$ be the 
$\bold Q$-subspace of $\bold C$ generated by $p_n(H)$ for $H \in \Cal A_S$.
Then it is easy to see that for finite subsets $S, T$ such that 
$T \subset S$, we have $p_n(T) \subset p_n(S)$.

\subsection{Periods and multiple zeta values}
In this section, we construct topological cycles of
relative homology group and compute the natural pairing
$$
H^n_{DR}(Y^0, j_!\bold Q) \otimes 
H_n^B(Y^0, \text{ mod }\bold B^0) \to \bold C,
$$
where $\bold B^0$ is the normal crossing divisor defined in \S
\ref{subsec:relative de Rham cohom}
and $j :Y^0 - \bold B^0 \to Y^0$ be the natural inclusion.
For an element 
$\gamma$ in $\pi_1(\bold A^1 -\{ 0,1 \},\frac{1}{2})$,
we define a path $\bar \gamma$ connecting $0$ and $1$
by $[\frac{1}{2},1]\circ \gamma \circ [0,\frac{1}{2}]$, where
$\circ$ denotes the composite of paths. It is well defined up to
homotopy equivalece. We define a continuous map $\delta_{\gamma}$ from 
$\Delta_n = \{ 0< t_1 < \cdots < t_n < 1\}$ to $Y^0$ by
$$
\delta_{\gamma}(t_1, \dots ,t_n) = 
(\bar\gamma (t_1), \dots , \bar\gamma (t_1)).
$$
Since the variety $(\bold A^1 -\{0,1\})^n$ can be
identified with an open set of $Y$, the map $\delta_{\gamma}$
can be lifted to a map $\delta_{\gamma, Y}$ to $Y$.
\begin{proposition}
The closure 
$\overline{Im (\delta_{\gamma, Y})}$ of the image 
$Im (\delta_{\gamma,Y})$ is contained in $Y^0$.
The boundary $\overline{Im (\delta_{\gamma, Y})} -
Im (\delta_{\gamma, Y})$ is contained in $\bold B^0$. 
\end{proposition}
\begin{proof}
For the first assertion, it is enough to prove that 
$\overline{Im (\delta_{\gamma, Y})}\cap E_i^{pr}$,
$\overline{Im (\delta_{\gamma, Y})}\cap F_i^{pr}$
($i=0, \dots ,n-2$) and $\overline{Im (\delta_{\gamma, Y})}\cap B_i^{pr}$
($i=0, \dots ,n$) does not intersect with $D_Y$.

(1) {\bf Proof of $\overline{Im (\delta_{\gamma, Y})}\cap E_i^{pr}
\cap D_Y=\emptyset$.}
By the natural map $\pi_i : E_i^{pr} \to Y(x_{n-i+1}, \dots ,x_n)$,
$\pi_i(\overline{Im (\delta_{\gamma, Y})}\cap E_i^{pr})$
does not intersect with 
\linebreak
$D_Y(x_{n-i+1}, \dots ,x_n)$
by the inductive hypothesis. For 
$y \in \pi_i(\overline{Im (\delta_{\gamma, Y})}\cap E_i^{pr})$,
$\pi_i^{-1}(y) \cap Im (\delta_{\gamma, Y}) \cap E_i^{pr,0}$
is equal to 
$$
0 \leq \eta_i \leq 1 \quad (i=1, \dots , n-i-1)
$$
with the coordinate of (\ref{fiber equation}).
Therefore $\overline{Im (\delta_{\gamma, Y})}\cap E_i^{pr}$ does 
not intersect with $D_Y$.

(2) 
{\bf Proof of $\overline{Im (\delta_{\gamma, Y})}\cap F_i^{pr}
\cap D_Y=\emptyset$.}
The proof is similar.

(3) 
{\bf Proof of $\overline{Im (\delta_{\gamma, Y})}\cap B_i^{pr}
\cap D_Y=\emptyset$.}
Since 
$B_i^{pr} \simeq$
\linebreak
$Y(x_1, \dots ,x_{i-1}, x_{i+1}, \dots ,x_n)$
for $i=1, \dots ,n$ and 
$B_0^{pr} \simeq Y(x_2, \dots ,x_n)$, 
\linebreak
$\overline{Im (\delta_{\gamma, Y})}\cap B_i^{pr}$ does not
intersect with $D_Y \cap B_i^{pr}$ by the inductive hypothesis.

The second assertion is a direct consequence of the first assertion.
\end{proof}
The relative cycle $\overline{Im (\delta_{\gamma, Y})}$ is denoted by
$\bar \delta_{\gamma}$
For an element $\eta =\sum_{\gamma}a_{\gamma}\gamma \in 
\bold Q[\pi_1(\bold A^1 -\{ 0,1 \},\frac{1}{2})]$, we define an element
$\bar\delta (\eta)$ by
$\sum_{\gamma}a_{\gamma}\bar\delta_{\gamma}$ in 
\linebreak
$H_n(Y^0,\text{ mod }\bold B^0)$. Then we have a linear map
$$
\bar\delta :\bold Q[\pi_1(\bold A^1 -\{ 0,1 \},\frac{1}{2})] \to
H_n(Y^0,\text{ mod }\bold B^0).
$$
Since the filtration induced by spectral sequence (\ref{eqn:spect}) coincides
with the weight spectral sequence, $F^n(H^n_{DR}(Y^0,j_!\bold Q_V))$
is generated by differential form 
$\omega (1, \dots ,n;\epsilon_1, \dots ,\epsilon_n)$
such that $\epsilon_i=0,1$ $i=2, \dots , n-1$ and $\epsilon_1=1, \epsilon_n=0$.
By the definition of $\bar\delta$, 
$$
(-1)^l\langle\bar\delta (1),\omega (1, \dots ,n;\epsilon_1, \dots ,\epsilon_n)\rangle=
\zeta (k_1, \dots ,k_l),
$$
where $\epsilon_j = 1$ for $j=h_0+1, h_1+1, \dots ,h_{l-1}+1$ 
and $\epsilon_j=0$ otherwise.
Here we use the correspondence between $k_i$ and $h_i$
in the introduction.
Therefore $\zeta (k_1, \dots ,k_l) \in p_n(H^n(j_!\bold Q_V, Y^0))$,
and we have the following proposition.
\begin{proposition}
\label{prop:real str and mzv}
$L_n \subset p_n(H^n(Y^0, j_!\bold Q_V)) \cap \bold R$.
\end{proposition}

\subsection{Periods of subquotients of $j_!\bold Q_V[n]$}
\label{subsec:period of subq}
In this subsection, we introduce an inductive structure
to compute the periods of subquotients of 
\linebreak
$H^n(Y^0, j_!\bold Q_V)$.

We define a divisor ${\bold B'}^0= \cup_{i=1}^n \tilde B_i^0$ 
and $j'$ as the
open immersion from $V'=Y^0-{\bold B'}^0$ to $Y^0$. 
The complement $\bold B^0 - {\bold B'}^0$ is denoted by $V_{n-1}$
and the open immersion $V_{n-1} \to \tilde B_0^0$ is denoted by $j_{n-1}$.
It is easy to see that $j'_!\bold Q_{V'}\mid_{\tilde B_0^0} 
\simeq j_{n-1,!}\bold Q_{V_{n-1}}$.
Therefore we have the following distinguished triangle
\begin{align}
\label{eqn:triang inductive}
j_{!}\bold Q_{V} \to
j'_{!}\bold Q_{V'} \to 
j_{n-1,!}\bold Q_{V_{n-1}}\overset{+1}\to .
\end{align}
$\tilde B_0^0$ is homotopy equivalent to 
$B_0^{pr,0} \simeq Y(x_2, \dots ,x_n)^0$, and $V_{n-1}$ is equal
to $\tilde B_0^0-\cup_{i=1}^n \tilde B_i^0$.
\begin{proposition}
\label{prop:induct}
The morphism
\begin{equation}
\label{eqn:inductive motives}
j_{n-1,!}\bold Q_{V_{n-1}}[n-1] \to j_{!}\bold Q_{V}[n]
\end{equation}
arising from triangle (\ref{eqn:triang inductive}) 
is injective and the image is
identified with 
\linebreak
$W_{n-1}j_{!}\bold Q_{V}[n]$.
\end{proposition}
\begin{proof}
Let $\oplus_{\# S=n+1-*, 0 \in S\subset [0,S]}\bold Q_{\tilde B_S^0}$ and
$\oplus_{\# S=n+1-*, S\subset [1,S]}\bold Q_{\tilde B_S^0}$ 
be the complices defined
similarly to (\ref{eqn:complex}).
Then the triangle 
$$
j_{n-1,!}\bold Q_{V_{n-1}}[-1] \to
j_!\bold Q_V \to j'_!\bold Q_{V'} 
$$
is isomorphic to
\begin{align*}
\oplus_{\# S=n+1-*, S\subset [1,S]}\bold Q_{\tilde B_S^0}
& \to
\oplus_{\# S=n+1-*, S\subset [0,S]}\bold Q_{\tilde B_S^0} \\
& \to
\oplus_{\# S=n+1-*, 0 \in S\subset [0,S]}\bold Q_{\tilde B_S^0}.
\end{align*}
By taking the Hodge realization of this triangle, morphism
(\ref{eqn:inductive motives}) is an injective morphism
between abelian objects.
\end{proof}
This proposition gives an inductive structre for relative cohomologies.
We claim the comptibility of this inductive structure and 
the homomorphism $\bar\delta$.
\begin{proposition}
\label{prop:dual induct}
Let $H_{n}(Y^0, \text{ mod }\bold B^0) \to 
H_{n-1}(\tilde B_0^0, \text{ mod }{\bold B'}^0 \cap \tilde B_0^0)$
be the dual of the Hodge realization of
the homomorphism given in Proposition \ref{prop:induct}.
Then the following diagram commutes.

\begin{tabular}{ccc}
$\bold Q [\pi_1(\bold A^1 -\{ 0,1 \},\frac{1}{2})]$ & $\to$ &
$\bold Q [\pi_1(\bold A^1 -\{ 0,1 \},\frac{1}{2})]$  \\
$\bar\delta \downarrow$  & & $\bar\delta \downarrow$ \\
$H_{n}(Y^0, \text{ mod }\bold B^0)$  & $\to$ &
$H_{n-1}(\tilde B_0^0, \text{ mod }{\bold B'}^0 \cap \tilde B_0^0)$. \\
\end{tabular}

Moreover, we have
$$
H_{n-1}(\tilde B_0^0, \text{ mod }{\bold B'}^0 \cap \tilde B_0^0)
\simeq
H_{n-1}(B_0^{pr,0}, \text{ mod }\bold B_{n-1}^0 ),
$$
where $\bold B^0_{n-1}$ is defined by
$$
B^{pr,0}_0 \cap ( E_0^{pr,0} \cup \cdots \cup E_{n-2}^{pr,0} \cup
B_2^{pr,0} \cup \cdots \cup B_{n-1}^{pr, 0} \cup \tilde B_n^{0}).
$$

\end{proposition}
We compute the period of subquotients of $H^n(Y^0, j_!\bold Q_V)$ via
formulas for iterated integrals.
Let $\omega_1, \dots , \omega_k$ be 1-forms on $[0,1]$.
The iterated integral
$\displaystyle\int_0^x\omega_1 \cdots  \omega_k$
is inductively defined by
$$
\int_0^x\omega_1 \cdots  \omega_k =
\int_0^x[\omega_1(u) 
(\int_0^u\omega_2 \cdots  \omega_k )].
$$
For the properties of iterated integrals, see \cite{Ch}.
Let $\gamma$ be a path from $[0,1]$ to $\bold A^1-\{ 0,1\}$ and
$\omega_1, \dots , \omega_k$ be holomorphic 1-form on $\bold A^1-\{ 0,1\}$.
The iterated integral along the path $\gamma$ is defined by
$$ 
\int_{\gamma}\omega_1 \dots  \omega_k
=\int_0^1\gamma^{*}(\omega_1) \dots  \gamma^{*}(\omega_k).
$$
It is known that it depends only on the homotopy quivalece class of
$\gamma$.
For a formal $\bold Q$-linear combination $\gamma = \sum_i c_i\gamma_i$
of paths from $p_0$ to $p_1$,
we define $\displaystyle\int_{\gamma} \omega_1\cdots \omega_k$
by $\displaystyle\sum_i c_i \int_{\gamma_i} \omega_1\cdots \omega_k$.

Let $\gamma$ (resp. $\delta$) be paths connecting $b$ and $c$ ($a$ and $b$)
and $\gamma \cdot \delta$ be the composite of these paths.
Then we have the following coproduct formula.
$$
\int_{\gamma \cdot \delta}\omega_1 \cdots \omega_k
=
\int_{\gamma}\omega_1 \cdots \omega_k
+ \sum_{i=1}^{k-1}(\int_{\gamma}\omega_1 \cdots \omega_i)
\cdot (\int_{\delta}\omega_{i+1} \cdots \omega_k)
+ \int_{
\delta}\omega_1 \cdots \omega_k.
$$

Let $\gamma$ be a path from $1$ to $0$ such that $\gamma=x$
if $x \in [0,\epsilon)$ or $(1-\epsilon, 1]$ for a sufficiently small 
$\epsilon$. Then
$$
\int_{\gamma}\frac{dx}{x-e_n}\cdots\frac{dx}{x-e_1}
=
\lim_{\delta\to 0}\int_{\gamma\mid [\delta, 1-\delta]}
\frac{dx}{x-e_n}\cdots\frac{dx}{x-e_1}
$$
exists if $e_1=1, e_n=0$. It is easy to see that
$$
\int_{[\frac{1}{2},1]\gamma [0,\frac{1}{2}]}\frac{dx}{x-e_n}\cdots\frac{dx}{x-e_1}
=
\langle\bar\delta (\gamma),
\omega (1, \dots ,n ;e_1, \dots ,e_n\rangle).
$$

Let $\rho_0$ (resp. $\rho_1$) be a small loop runinging around $0$ (resp. $1$)
with base point $\frac{1}{2}$.
The group ring 
$\bold Q [\pi_1(\bold A^1 -\{ 0,1 \},\frac{1}{2})]$ is denoted by $R$ 
and the augumentation ideal $Ker (R \to \bold Q)$ is denoted by $I$.
\begin{proposition}
\label{prop:iterate formula}
We put $\alpha=[\frac{1}{2},1]$, $\beta=[0,\frac{1}{2}]$. Suppose that
$e_i=0,1$ for $i=2, \dots, n-1$, $e_1=0$ and $e_n=1$.

\begin{enumerate}
\item
If $\gamma$ is contained in $I^{n+1}$, then
$$
\int_{\alpha \gamma \beta}
\frac{dx}{x-\epsilon_1}\cdots \frac{dx_n}{x-\epsilon_n}
=
0
$$
\item
$$
\int_{\alpha \cdot (g_1-1) \cdots (g_n-1) \cdot \beta}
\frac{dx}{x-\epsilon_1}\cdots \frac{dx_n}{x-\epsilon_n}
=
\begin{cases}
(2 \pi i)^n & \text{ if } g_i= \rho_{\epsilon_i} \text{ for all } i\\
0 & \text{ otherwise }\\
\end{cases}
$$
\item
Let $g_i \in \{ \rho_0, \rho_1\}$, $\epsilon_k = 0, 1$ 
and $\epsilon_1=1, \epsilon_n=0$.
Then we have
$$
\int_{\alpha \cdot (g_1-1) \cdots (g_{n-1}-1) \cdot \beta}
\frac{dx}{x-\epsilon_1}\cdots \frac{dx_n}{x-\epsilon_n}
\in \frac{(2\pi i)^n}{2}\bold Z.
$$
\end{enumerate} 
\end{proposition}
\begin{proof}
1. This equality comes from the coproduct formula.

2. By the coproduct formula, we have
$$
\int_{\alpha \cdot (g_1-1) \cdots (g_n-1) \cdot \beta}\omega_1\cdots \omega_n
=
\prod_{i=1}^n\int_{g_i}\omega_i,
$$
and the proposition follows from this equality.

3. 
We use the notation
$\displaystyle(\gamma,\epsilon_1\cdots \epsilon_k)=
\int_{\gamma}\frac{dx}{x-\epsilon_1}\cdots \frac{dx_k}{x-\epsilon_k}$.
By the coproduct formula, we have
\begin{align}
\label{eqn:subquot period}
& \int_{\alpha \cdot (g_1-1) \cdots (g_{n-1}-1) \cdot \beta}
\frac{dx}{x-\epsilon_1}\cdots \frac{dx_n}{x-\epsilon_n} \\
= &
(\alpha, \epsilon_1)
\prod_{i=1}^{n-1}(g_i,\epsilon_{i+1}) 
+
(\beta, \epsilon_{n})
\prod_{i=1}^{n-1}(g_i,\epsilon_{i})  \nonumber \\
& + \sum_{k=1}^{n-1}
\prod_{i=1}^{k-1}(g_i,\epsilon_{i})\cdot
(g_k,\epsilon_k \epsilon_{k+1})\cdot
\prod_{i=k+1}^{n-1}(g_i,\epsilon_{i+1}). \nonumber
\end{align}
Let $p$ (resp. $q$)
be the minimal number such that $g_{p+1}\neq \rho_{\epsilon_{p+1}}$
(resp. maximal number such that $g_{q-1} \neq \rho_{\epsilon_{q}}$).
If $p+1 \leq q-1$, i.e. there exists no $k$ 
such that $q-1 \leq k \leq p+1$,
then all the terms in (\ref{eqn:subquot period}) vanish.
Therefore we have the proposition. 
We may assume that $q-1 < p+1$.
We can easily see that $q>1$ or $p<n-1$.
Suppose that $q>1$. If $q-1 < k <p+1$, (resp. $q-1>k$ or $k > p+1$)
then
$\prod_{i=1}^{k-1}(g_i,\epsilon_{i})\cdot
(g_k,\epsilon_k \epsilon_{k+1})\cdot
\prod_{i=k+1}^{n-1}(g_i,\epsilon_{i+1})$ is equal to
$(2 \pi i)^n/2$ (resp. $0$).

(Case 1) If  $p= n-1$, we have
$$
\prod_{i=1}^{q-2}(g_i,\epsilon_{i})\cdot
(g_{q-1},\epsilon_{q-1} \epsilon_{q})\cdot
\prod_{i=q}^{n-1}(g_i,\epsilon_{i+1}) + 
(\beta, \epsilon_{n})
\prod_{i=1}^{n-1}(g_i,\epsilon_{i}) = 0.
$$
Therefore the sum ($\ref{eqn:subquot period}$)
is contained in $\displaystyle\frac{(2 \pi i)^n}{2} \bold Z$.

(Case 2) If $p < n-1$, we have
\begin{align*}
& \prod_{i=1}^{q-2}(g_i,\epsilon_{i})\cdot
(g_{q-1},\epsilon_{q-1} \epsilon_{q})\cdot
\prod_{i=q}^{n-1}(g_i,\epsilon_{i+1}) \\
& + 
\prod_{i=1}^{p}(g_i,\epsilon_{i})\cdot
(g_{p+1},\epsilon_{p+1} \epsilon_{p+2})\cdot
\prod_{i=p+2}^{n-1}(g_i,\epsilon_{i+1})= 0.
\end{align*}
Again the sum ($\ref{eqn:subquot period}$) is 
contained in $\displaystyle\frac{(2 \pi i)^n}{2} \bold Z$.

For the case $q=1$ and $p< n-1$, the proof is similar.
\end{proof}

Set $h_i = \rho_i -1 \in R$.
\begin{corollary}
The restriction of the morphism $\bar\delta$ 
in Proposition \ref{prop:dual induct}
to $\bold Q\cdot 1 \oplus h_0 R h_1$
is surjective and $H_{n}(Y^0, \text{ mod }\bold B^0)$
is isomorphic to 
$\bold Q\cdot 1 \oplus$ 
\linebreak
$h_0 R h_1
/h_0 I^{n-1} h_1$ via this homomorphism.
Moreover the weight filtration coincides with that induced
from the power of $I$.
\end{corollary}
\begin{proof}
We prove the corollary by induction.
For $n=2$, we can prove the proposition directly.
By the duality and the first statement of 
Proposition \ref{prop:iterate formula}, we have
$\bar\delta (h_0 I^{n-1}h_1)=0$.
The following diagram is commutative by Proposition \ref{prop:dual induct}.

\begin{tabular}{ccc}
0  & & 0
\\
$\downarrow$ & & $\downarrow$
\\
$h_0 I^{n-2}h_1/ h_0I^{n-1}h_1$
& $\overset{(1)}\to$ &
$W_{-n}H_n(Y^0, \text{ mod }\bold B^0)$
\\
$\downarrow$ & & $\downarrow$
\\
$\bold Q\cdot 1 \oplus h_0 R h_1/ h_0 I^{n-1} h_1$
& $\overset{(2)}\to$ &
$H_{n}(Y^0, \text{ mod }\bold B^0)$
\\
$\downarrow$ & & $\downarrow$
\\
$\bold Q\cdot 1 \oplus h_0 R h_1/ h_0 I^{n-2} h_1$
& $\overset{(3)}\to$ &
$H_{n-1}(B_0^0,\text{ mod }{\bold B'}^0\cap B_0^0)$ \\
$\downarrow$ & & $\downarrow$
\\
0 & & 0
\end{tabular}

Here the columns are exact.
Homomorphism $(3)$ is an isomorphism by the assumption of induction.
Using the duality, homomorphism $(1)$ is an isomorphism by
Proposition \ref{prop:iterate formula}.2.
\end{proof}
\begin{theorem}
\label{thm:level one split}
Let $H^n=H^n(Y^0, j_{!}\bold Q_{V})$. Then the exact sequence
\begin{align}
\label{eqn:lev1 geom}
0 \to W_{k-1}(H^n)/W_{k-2}(H^n) 
& \to W_k(H^n)/W_{k-2}(H^n) \\
& \to W_k(H^n)/W_{k-1}(H^n) \to 0 \nonumber
\end{align}
splits as mixed Hodge structures for all $k$.
\end{theorem}
\begin{proof}
By the inductive structure introduced in 
Proposition \ref{prop:dual induct},
we may assume $k=n$.
We compute the extension class of (\ref{eqn:lev1 geom}) in 
\linebreak
$Ext^1(Gr_n^W(H_n),Gr_{n-1}^W(H_n))$ by the recipe 
in \S \ref{sebsec:comparison extension grp}.
By Proposition \ref{prop:iterate formula}.3,
the differential form 
$\frac{1}{(2 \pi i)^n}\omega(1,\dots ,n;1,\epsilon_2 ,\dots, \epsilon_{n-1},0)$
with $\epsilon_i=0,1$ ($i=2, \dots , n-1$) forms a $\bold Q$ base of
$W_n(H^n)/W_{n-1}(H^n)$ via the isomorphism
$$
F^n(H^n\otimes \bold C) \simeq (W_n(H^n)/W_{n-1}(H^n))\otimes \bold C.
$$
On the other hand, $\pi (\bar\delta (h_0(g_2-1)\cdots (g_{n-2}-1)h_1 ))$
form a basis in
\linebreak
$H_{n-1}(B_0^{pr,0},\text{ mod }\bold B_{n-1}^0)$
under the morphism 
$$
\pi :
H_{n}(Y^0, \text{ mod }\bold B^0)
\to H_{n-1}(B_0^{pr,0},\text{ mod }\bold B_{n-1}^0).
$$
Therefore the pairing 
\begin{align*}
& \langle\bar\delta (h_0(g_2-1)\cdots (g_{n-2}-1)h_1 ),
\frac{1}{(2 \pi i)^n}
\omega(1,\dots ,n;1,\epsilon_2 ,\dots, \epsilon_{n-1},0)\rangle \\
= &
\frac{1}{(2 \pi i)^n}\int_{h_0(g_2-1)\cdots (g_{n-2}-1)h_1}
\frac{dx}{x}\frac{dx}{x-e_{n-1}}\cdots \frac{dx}{x-e_{2}}\frac{dx}{x-1}
\in \frac{1}{2}\bold Z
\end{align*}
and the extension class of (\ref{eqn:lev1 geom}) vanishes.

\end{proof}
We define $\Cal N_n^{\geq 2}$ by
$$
\Cal N_n^{\geq 2} = \{ \{a_1, \dots , a_n\} \mid
0 \leq a_1 < \cdots < a_n \leq n
\text{ and } a_{i+1}-a_i \geq 2 \text{ for all } i
\}.
$$
\begin{corollary}
\label{cor:mzv split in level 1}
Let $H^n=j_!\bold Q_V[n]\in \Cal A_{TM}$.
Exact sequence (\ref{eqn:lev1 geom}) splits in the category 
$\Cal A_{TM}$.
$j_!\bold Q_V[n]$ is a subquotient of a direct product of objects
in $\Cal A_S$ with $S \in \Cal N_n^{\geq 2}$. 
\end{corollary}
\begin{proof}
The first statement is a direct consequence of Theorem \ref{thm:level one split}
and Proposition \ref{prop:compatiblity ext one}.
We prove the second statement by the induction on $n$.
Since $H^n=j_!\bold Q_V[n]$ splits in level $(n-1,n)$, 
by Proposition \ref{prop: decomp},
there exist
an object $B_1$ and $B_2$ in $\Cal A$ such that the following
statement holds.
(1) $H^n$ is a subquotient of
$B_1 \oplus B_2$.  (2) $B_1$ is an object of 
$\Cal A_{[0,n-2]\cup \{ n\}}$.
$W_{n-2}B_1$ is isomorphic to 
the direct sum of copies of $W_{n-2}(H^n)$.
(3) $W_{n-1}B_2 = B_2$ is isomorphic to 
the direct sum of copies of $W_{n-1}(H^n)$.
By the assumption of induction, $W_{n-2}(H^n)$ (resp. $W_{n-1}(H^n)$)
is a subquotient of a direct sum of
objects in $A_S \in \Cal A_S$ with 
$S \in \Cal N_{n-2}^{\geq 2}$ (resp. $S \in \Cal N_{n-1}^{\geq 2}$).
Therefore $B_1$ and $B_2$ are subquotients of a direct sum of objects
in $A_S \in \Cal A_S$ with 
$S \in \Cal N_{n}^{\geq 2}$ and $S \in \Cal N_{n-1}^{\geq 2}$, respectively.
\end{proof}

\section{Proof of Main Theorem}

We use the same notations $\Cal N_n$ and $\Cal N_n^{\geq 2}$
as in \S \ref{subsec:generator} and \S \ref{subsec:period of subq}
respectively.
\begin{lemma}
\label{lem:reduction odd}
Let $S=\{a_1< \cdots < a_l\}$. 
\begin{enumerate}
\item
If $S \in \Cal N_n$,
then 
$\dim p_n(S)/\sum_{T\subsetne S}p_n(T) \leq 1$.
\item
If $S \in \Cal N_n^{\geq 2}-\Cal N_n$,
then $\dim p_n(S)/\sum_{T\subsetne S}p_n(T)=0$.
\end{enumerate}
\end{lemma}
\begin{proof}
If every $a_{k+1}-a_{k}$ is an odd number greater than 1,
$\Cal A_S =<M_S,\Cal T>_{T\subsetne S}$. 
Since $M_S/W_{a_1}M_S \in \Cal A_{S-\{ a_1 \}}$,
the set $\{ \langle\gamma, \omega\rangle\mid
\omega \in F_{n}M_S, \gamma \in W_{-a_2}M_S^*\}$
is contained in $p_n(S-\{ a_1 \})$.
Since $\dim W_{-a_1}M_S^*/W_{-a_2}M_S^* =1$ and 
\linebreak
$\dim F^n(M_S) \leq 1$,
$p_n(S)/\sum_{T\subsetne S}p_n(T)$ is at most one dimensional.

Suppose that $a_{i+1}-a_i$ is even. Let $A\in \Cal A_S$.
Since
\linebreak
$Ext^1(\bold Q(-a_i),\bold Q(-a_{i+1}))=0$, $A$ splits in level $(a_i, a_{i+1})$.
Therefore $A$ can be written as a subquotient of a direct sum $B_1 \oplus B_2$
with $B_1 \in \Cal A_{S-\{ a_i\}}$ and $B_2 \in \Cal A_{S-\{ a_{i+1}\}}$
by Proposition \ref{prop: decomp}. Therefore we have the second statement.
\end{proof}
By Corollary \ref{cor:mzv split in level 1},
the space $p_n(H^n(j_!\bold Q_V))$ is contained in the space
\linebreak
$\sum_{S \in \Cal N_n^{\geq 2}}p_n(S)$.
Since $\sum_{S \in \Cal N_n^{\geq 2}}p_n(S)=\sum_{S \in \Cal N_n}p_n(S)$
by Lemma \ref{lem:reduction odd}, 
\linebreak
$p_n(H^n(j_!\bold Q_V))$
is contained in $\sum_{S \in \Cal N_n}p_n(S)$.

Let $\op (a)$ be the cardinality of the set 
$$
\{ (b_1, b_2, \cdots) \text{ (ordered) }\mid b_i 
\text{ is an odd integer greater than 1 },\sum_ib_i = a \}.
$$
Since $\# \{S \in \Cal N_n, S \ni a, n, S\subset [a,n] \}=\op (n-a)$,
we have
$$
\dim (\sum_{\substack{S \in \Cal N_n \\  S \subset [a,n]}}p_n(S)
/\sum_{\substack{S \in \Cal N_{n} \\ S \subset [a+1,n]}}p_n(S) )
\leq \op (n-a)
$$
by Lemma \ref{lem:reduction odd}.
Moreover the complex conjugate $c$ acts on the space \linebreak
$\sum_{S \in \Cal N_n, S \subset [a,n]}p_n(S)
/\sum_{S \in \Cal N_{n},S \subset [a+1,n]}p_n(S)
$ by the $(-1)^a$-multiplication.
By Lemma \ref{prop:real str and mzv},
$L_n$ is a subset of $p_n(H^n(j_!\bold Q_V))$
and invariant under the complex conjugation.
We have
\begin{align*}
 \dim L_n & \leq p_n(H^n(j_!\bold Q_V)) \cap \bold R \\
& =  
\sum_{\substack{0 \leq a \leq n \\ a\equiv 0 (\text{ mod } 2)}} 
\dim (\sum_{\substack{S \in \Cal N_n \\  S \subset [a,n]}}p_n(S)
/\sum_{\substack{S \in \Cal N_{n-1} \\ S \subset [a+1,n]}}p_n(S) )  \\
& \leq 
\sum_{\substack{0 \leq a \leq n \\ a\equiv 0 (\text{ mod } 2)}} 
\op (n-a). 
\end{align*}
Therefore it is enough to prove the following lemma.
\begin{lemma}
Let $d_n$ be the number defined in the introduction. Then we have
$$
d_n =\sum_{\substack{0 \leq a \leq n \\ a\equiv 0 (\text{ mod } 2)}} 
\op (n-a).
$$
\end{lemma}
\begin{proof}
We consider variables $u_3, u_5, \dots $, where the degree of $u_i$
is $i$. We put $V=\oplus_{i=1}^{\infty}\bold Qu_{2i+1}$
and $\bold V = \oplus_{i=0}^{\infty} V^{\otimes i}$
and $U = \oplus_{j=0}^{\infty}\bold V \cdot w^{i}$, where
$w$ is a variable of degree 2. Then the poincare series of $U$
is equal to
$$
\sum_{k=0}^{\infty}[\sum_{j=0}^{\infty}(\sum_{i=1}^{\infty}t^{2i+1})^j]t^{2k}
=\frac{1}{1-t^2-t^3}.
$$
Therefore the coefficient $d_n$ is equal to 
$\displaystyle\sum_{\substack{a \geq 0 \\ a\equiv 0 (\text{ mod } 2)}} 
\op (n-a)$.
\end{proof}

\end{document}